\RequirePackage{fix-cm}
\documentclass[smallextended, envcountsame]{svjour3}
\smartqed
\journalname{}

\usepackage[unicode = false,
    pdftoolbar = true,
    colorlinks = true,
    linkcolor = blue,
    citecolor = blue,
    filecolor = black,
    urlcolor = blue,
    breaklinks = true]{hyperref}

\usepackage[utf8]{inputenc}
\usepackage{listings}           
\usepackage{enumitem}           
\usepackage{amsmath}            
\usepackage{amssymb}            
\usepackage{xcolor}             
\usepackage{graphicx}           
\usepackage[ruled]{algorithm2e}
\usepackage{verbatim}
\usepackage{geometry}
\usepackage[justification=centering]{caption}

\definecolor{gab}{HTML}{50c878}
\definecolor{pat}{HTML}{bf68ff}
\definecolor{sham}{HTML}{7faaff}

\allowdisplaybreaks[4]

\title{Error Analysis of Surrogate Models Constructed through Operations on Sub-models}
\author{Yiwen Chen Warren Hare Gabriel Jarry-Bolduc }
\institute{Yiwen Chen \at
Department of Mathematics, University of British Columbia, Kelowna,\\
              British Columbia, Canada.  Chen's research is partially funded by the MITACS Globalink program and the Natural Sciences and Engineering Research Council (NSERC) of Canada, Discover Grant \#2018-03865.\\
\email{yiwen.chen@ubc.ca}
\and
Warren Hare \at
 Department of Mathematics, University of British Columbia, Kelowna, \\
              British Columbia, Canada.  Hare's research is partially funded by the Natural Sciences and Engineering Research Council (NSERC) of Canada, Discover Grant \#2018-03865.\\
              \email{warren.hare@ubc.ca}   
\and
Gabriel Jarry-Bolduc\at
           Department of Mathematics, University of British Columbia, Kelowna, \\
              British Columbia, Canada.  Jarry-Bolduc's research is partially funded through the Natural Sciences and Engineering Research Council (NSERC) of Canada, Discover Grant \#2018-03865.\\
              \email{gabjarry@alumni.ubc.ca} 
              }
\date{\today}

\begin{document}

    \maketitle
    
    \begin{abstract}
    Model-based methods are popular in derivative-free optimization (DFO). In most of them, a single model function is built to approximate the objective function. This is generally based on the assumption that the objective function is one blackbox. However, some real-life and theoretical problems show that the objective function may consist of several blackboxes. In those problems, the information provided by each blackbox may not be equal. In this situation, one could build multiple sub-models that are then combined to become a final model. In this paper, we analyze the relation between the accuracy of those sub-models and the model constructed through their operations. We develop a broad framework that can be used as a theoretical tool in model error analysis and future research in DFO algorithms design. 
    \end{abstract}

    \section{Introduction}

    Derivative-Free Optimization (DFO) is the mathematical study of algorithms for continuous optimization that do not use first-order information \cite{audet2017derivative}. In general, DFO methods can be categorized into model-based methods and direct search methods. In model-based methods, a model function is built to approximate the objective function. If the model is ‘well-built’, then it can behave well in many classic optimization algorithms\cite{audet2020model}. Therefore, methods to construct models are of great interest in DFO.
    
    Most model-construction methods assume that the objective function is a single blackbox. However, in some real-life problems, the objective function contains several blackboxes. For example, Rashki \cite{rashki2018hybrid} presents an optimization problem in structural reliability. Rashki first shows that the blackbox objective function can be written as the product of several other blackbox functions. Then, he builds models for each of them separately, and demonstrates that the restructured problem provides the objective function value at reduced computational cost. 
    
    Similar approaches are also used in other fields. In circuit design, solving fractional-order differential equations in the time domain is mathematically demanding \cite{hwang2002note,wu2004new}. In \cite{maundy2011practical}, the authors rewrite the objective function as the product of two functions according to frequencies and approximate each of them separately using continuous fraction expansions \cite{chen2004continued,chen2002discretization}. In \cite{kandasamy2013simulation}, the authors aim to maximize the overall waterjet-hull system efficiency by changing shape parameters of the ship model. They decompose the overall efficiency into several partial efficiency and evaluate independently by experiments and measurements. In the above problems, the different blackboxes do not have the same amount of information and models are built for each component blackbox separately. 
    
    Similarly, some theoretical work is also based on the assumption that the component functions involved are not equally expensive. Khan et al. \cite{khan2018manifold} develop an algorithm for minimizing $F=\phi + h\circ f$ where $\phi$ is smooth with known derivatives, $h$ is known, nonsmooth piecewise linear function and $f$ is smooth but expensive to evaluate. Later in \cite{larson2020manifold}, Larson et al. aim to minimize $F=h\circ f$ where $h$ is nonsmooth but inexpensive to compute and $f$ is smooth but its Jacobian is not available. 
    
    In this paper, we consider the relation between the accuracy of sub-models and that of their operations (product, quotient, composition), which can be regarded as a model of objective function. Some previous work has focused on similar topics. In \cite{hare2020calculus,regis2015calculus}, the authors introduce calculus rules (product, quotient, chain) for simplex gradients. Based on these rules, gradient approximations techniques are developed \cite{hare2020discussion,hare2020error}. The authors of \cite{hare2020hessian} present a Hessian approximation technique and develop another Hessian approximation technique based on calculus rules. Hare \cite{hare2017compositions} examine the error bounds for the composition of a convex lower semi-continuous function with a smooth vector-valued function under the assumption that it is possible to provide 
   “fully linear models” for each component of the vector-valued function. 
    
    The goal of this paper is to develop a theoretical framework that can be used in model error analysis. All the previous work mentioned above are encompassed and can be regarded as a specific case in our framework. However, it should be noted that our framework does not include accuracy analysis of sub-gradients. Therefore, the error bounds for cases including nonsmooth component functions \cite{hare2013derivative,larson2016manifold} are not covered.
    
    This paper is structured as follows. In section 2, we introduce some basic definitions and give out some important lemmas that are used throughout the paper. In sections 3-5, we provide the relation between the accuracy of sub-models and that of their product, quotient, and composition, respectively. Section 6 concludes our work and presents some possible follow-up research directions that may be profitable.

    \section{Definitions and Lemmas}
	
	In DFO, a model is built to perform as a surrogate of the objective function. If a model is ‘well-built’, then it can behave well in many classic optimization algorithms. In order to measure the quality of models, we introduce the order-$N$ terminology from \cite{hare2020discussion}. Notice that in all the definitions in this paper, $B_\Delta(x_0)$ is the closed ball centered at $x_0$ with radius $\Delta$.
	
	\begin{definition}{(Order-$N^f$ function accuracy)}\label{def:facc}
		Given $f\in\mathcal{C}^0, x_0\in\mathrm{dom}(f),\bar{\Delta}>0$, and $N^f\ge0$. We say that $\{\widetilde{f}_\Delta\}_{(0,\bar{\Delta}]}$ is a class of models of $f$ parameterized by $\Delta$ that provides order-$N^f$ function
		accuracy at $x_0$ if there exists a scalar $\kappa^f(x_0) \ge 0$ such that, given any $\Delta\in(0,\bar{\Delta}]$ the model $\widetilde{f}_\Delta$ satisfies
		\begin{align*}
			||f(x_0) - \widetilde{f}_\Delta(x_0)|| \le \kappa^f(x_0)\Delta^{N^f}.
		\end{align*}
		We say that $\{\widetilde{f}_\Delta\}_{(0,\bar{\Delta}]}$ is a class of models of $f$ parameterized by $\Delta$ that provides order-$N^f$ function
		accuracy near $x_0$ if there exists a scalar $\kappa^f(x_0) \ge 0$ such that, given any $\Delta\in(0,\bar{\Delta}]$ the model $\widetilde{f}_\Delta$ satisfies
		\begin{align*}
			||f(x) - \widetilde{f}_\Delta(x)|| \le \kappa^f(x_0)\Delta^{N^f}\ \ \ for\ all\ x \in B_\Delta(x_0).
		\end{align*}
	\end{definition}

	\begin{example}
	    Let $\{\widetilde{f}_\Delta\}_{(0,\bar{\Delta}]}$ be a class of models of $f$ satisfies $\widetilde{f}_\Delta(x_0)=f(x_0)$ for all $\Delta\in(0,\bar{\Delta}]$. Then $\{\widetilde{f}_\Delta\}_{(0,\bar{\Delta}]}$ provides order-$\infty$ function accuracy at $x_0$. Suppose all $\widetilde{f}_\Delta$ and $f$ are Lipschitz continuous. Suppose the constant $L \ge 0$ is such that for all $x,y\in\mathbb{R}^d$ and $\Delta\in (0,\bar{\Delta}]$ we have $||f(x)-f(y)||\le L||x-y||$ and $||\widetilde{f}_\Delta(x)-\widetilde{f}_\Delta(y)||\le L||x-y||$. Then $\{\widetilde{f}_\Delta\}_{(0,\bar{\Delta}]}$ provides order-1 function accuracy of $f$ near $x_0$.
	\end{example}
	
	\begin{proof}
	    For all $x\in B_{\Delta}(x_0)$ we have
	    \begin{align*}
	        ||f(x)-\widetilde{f}_\Delta(x)|| &= ||f(x)-f(x_0)+\widetilde{f}_\Delta(x_0)-\widetilde{f}_\Delta(x)||\\
	        &\le ||f(x)-f(x_0)|| + ||\widetilde{f}_\Delta(x_0)-\widetilde{f}_\Delta(x)||\\
	        &\le 2L\Delta.
	    \end{align*}
	    That is, $\{\widetilde{f}_\Delta\}_{(0,\bar{\Delta}]}$ provides order-1 function accuracy near $x_0$.
	    
	    $\hfill\qed$
	\end{proof}

    Analogously, we can define order-$N^g$ gradient accuracy and order-$N^h$ Hessian accuracy. Specifically, if $f\in\mathcal{C}^1$, then replacing all $f(\cdot)$ by $\nabla f(\cdot)$ we get order-$N^g$ gradient accuracy; if $f\in\mathcal{C}^2$, then replacing all $f(\cdot)$ by $\nabla^2 f(\cdot)$ we get order-$N^h$ Hessian accuracy.
    
    Note that a similar concept is presented by Conn et al \cite{conn20081geometry,conn20082geometry,conn2009introduction} where they investigated the framework of fully linear and fully quadratic models. These can be easily expressed in terms of order-$N$ terminology. A fully linear class of models provides order-2 function and order-1 gradient accuracy. A fully quadratic class of models provides order-3 function, order-2 gradient, and order-1 Hessian accuracy.
	
	Another definition that is important throughout this paper is derived from fundamental mathematical analysis. It describes the uniform boundedness of models.
	
	\begin{definition}{(Uniformly bounded models)}\label{def:ub}
		Given $f \in \mathcal{C}^0, x_0 \in \mathrm{dom}(f)$, and $\bar{\Delta} > 0$. We say that $\{\widetilde{f}_\Delta\}_{(0,\bar{\Delta}]}$ is a class of models of $f$ parameterized by $\Delta$ that is uniformly bounded at $x_0$ if there exists a scalar $M(x_0) \ge 0$ such that, given any $\Delta \in (0,\bar{\Delta}]$ the model $\widetilde{f}_\Delta$ satisfies
		\begin{align*}
			||\widetilde{f}_\Delta(x_0)|| \le M(x_0).
		\end{align*}
		We say that $\{\widetilde{f}_\Delta\}_{(0,\bar{\Delta}]}$ is a class of models of $f$ parameterized by $\Delta$ that is uniformly bounded near $x_0$ if there exists a scalar $M(x_0) \ge 0$ such that, given any $\Delta \in (0,\bar{\Delta}]$ the model $\widetilde{f}_\Delta$ satisfies
		\begin{align*}
			||\widetilde{f}_\Delta(x)|| \le M(x_0)\ \ \ for\ all\ x \in B_\Delta(x_0).
		\end{align*}
	\end{definition}

	Also from fundamental mathematical analysis we have that if $f\in\mathcal{C}^k$ where $k\in\{0,1,...\}$, then $f$ and its derivatives of order less or equal than $k$ are bounded in $B_{\bar{\Delta}}(x_0)$. When needed, we use $M_f(x_0)$, and $M_{\nabla^i f}(x_0)$ where $i\in\{1,...,k\}$ to denote these bounds.
    
    In this paper, unless stated otherwise, we assume that for all $i\in\{1,2,...\}$, $f_i$ is a function from $\mathbb{R}^d$ to $\mathbb{R}$, and $\{\widetilde{f}_{i,\Delta}\}_{(0,\bar{\Delta}]}$ is a class of models of $f_i$ parameterized by $\Delta$.
	
	Examining Definition \ref{def:facc}, we can easily see that higher accuracy orders indicate more accurate models. Therefore, it is clear that the following lemma holds.

	\begin{lemma}
		If the class of models $\{\widetilde{f}_\Delta\}_{(0,\bar{\Delta}]}$ provides order-$N$ accuracy of $f$, then for all $0\le N'\le N$, $\{\widetilde{f}_\Delta\}_{(0,\bar{\Delta}]}$ provides order-$N'$ accuracy of $f$.
	\end{lemma}

	Notice that the definition of uniformly bounded models is equivalent to order-$0$ function accuracy. Therefore, we immediately get the next corollary.
	
	\begin{corollary}
		If $f\in\mathcal{C}^0$ and $\{\widetilde{f}_\Delta\}_{(0,\bar{\Delta}]}$ provides order-$N^f$ function accuracy of $f$ at/near $x_0$, then $\{||\widetilde{f}_\Delta||\}_{(0,\bar{\Delta}]}$ is uniformly bounded at/near $x_0$. If $f\in\mathcal{C}^1$ and $\{\widetilde{f}_\Delta\}_{(0,\bar{\Delta}]}$ provides order-$N^g$ gradient accuracy of $f$ at/near $x_0$, then $\{||\nabla\widetilde{f}_\Delta||\}_{(0,\bar{\Delta}]}$ is uniformly bounded at/near $x_0$. If $f\in\mathcal{C}^2$ and $\{\widetilde{f}_\Delta\}_{(0,\bar{\Delta}]}$ provides order-$N^h$ Hessian accuracy of $f$ at/near $x_0$, then $\{||\nabla^2\widetilde{f}_\Delta||\}_{(0,\bar{\Delta}]}$ is uniformly bounded at/near $x_0$.
	\end{corollary}

    For the rest of the paper, we denote the bounds in the corollary above by $M_{\widetilde{f}}(x_0)$, $M_{\nabla\widetilde{f}}(x_0)$, and $M_{\nabla^2\widetilde{f}}(x_0)$, respectively. Remark, it can be shown that $M_{\widetilde{f}}(x_0) = \bar{\Delta}^{N^f}\kappa^f(x_0) + M_{f_1}(x_0)$, $M_{\nabla\widetilde{f}}(x_0) = \bar{\Delta}^{N^g}\kappa^g(x_0) + M_{\nabla f}(x_0)$, and $M_{\nabla^2\widetilde{f}}(x_0) = \bar{\Delta}^{N^h}\kappa^h(x_0) + M_{\nabla^2 f}(x_0)$.
    
    In the rest of the paper, we consider cases where the objective function can be expressed as the product, quotient, or composition of several blackbox functions, and provide error analysis for each case respectively.

    \section{Error Analysis of the Product of Model Functions}
	
	In this section, we consider the case where the objective function $F$ can be expressed as the product of several component functions. We first consider the 2-function case, then we extend the results to the $n$-function case.
	
	\subsection{Error Analysis of the Product of 2 Functions}
	
	Throughout this subsection we let $F=f_1f_2$ and $\widetilde{F}_\Delta=\widetilde{f}_{1,\Delta}\widetilde{f}_{2,\Delta}$. In the next theorem we analyze the function accuracy of $\{\widetilde{F}_\Delta\}_{(0,\bar{\Delta}]}$.
	
	\begin{theorem}{(Function accuracy at/near $x_0$)}\label{thm:2f}
		Let $f_1,f_2$ be in $\mathcal{C}^0$. Suppose for all $i\in\{1,2\}$, $\{\widetilde{f}_{i,\Delta}\}_{(0,\bar{\Delta}]}$ provides order-$N_i^f$ function accuracy of $f_i$ at/near $x_0$ with constant $\kappa_i^f(x_0)$. Let $N_F^f = \min\{N_1^f,N_2^f\}$. Then
		
		$\ $(i) if the accuracy provided is at $x_0$, then
		\begin{align*}
			|F(x_0) - \widetilde{F}_\Delta(x_0)| \le K_F^f(x_0)\Delta^{N_F^f};
		\end{align*}
	
		(ii) if the accuracy provided is near $x_0$, then for all $x\in B_\Delta(x_0)$
		\begin{align*}
			|F(x) - \widetilde{F}_\Delta(x)| \le K_F^f(x_0)\Delta^{N_F^f};
		\end{align*}
		where
		\begin{align*}
			K_F^f(x_0) = \min\{&M_{\widetilde{f}_2}(x_0)\kappa_1^f(x_0)\bar{\Delta}^{N_1^f-N_F^f} + M_{f_1}(x_0)\kappa_2^f(x_0)\bar{\Delta}^{N_2^f-N_F^f},\\
			&M_{f_2}(x_0)\kappa_1^f(x_0)\bar{\Delta}^{N_1^f-N_F^f} + M_{\widetilde{f}_1}(x_0)\kappa_2^f(x_0)\bar{\Delta}^{N_2^f-N_F^f}\}.
		\end{align*}
		In particular, $\{\widetilde{F}_\Delta\}_{(0,\bar{\Delta}]}$ provides order-$N_F^f$ function accuracy of $F$ at/near $x_0$.
	\end{theorem}
	
	\begin{proof}
		We only prove case (ii), since the proof of case (i) can easily be obtained by changing all $x$ to $x_0$. First notice
		\begin{align*}
			|F(x) - \widetilde{F}_\Delta(x)| =& |(f_1f_2 - \widetilde{f}_{1,\Delta}\widetilde{f}_{2,\Delta})(x)|\\
			=& |(f_1f_2 - f_1\widetilde{f}_{2,\Delta} + f_1\widetilde{f}_{2,\Delta} - \widetilde{f}_{1,\Delta}\widetilde{f}_{2,\Delta})(x)|\\
			=& |(f_1\widetilde{f}_{2,\Delta} - \widetilde{f}_{1,\Delta}\widetilde{f}_{2,\Delta} + f_1f_2 - f_1\widetilde{f}_{2,\Delta})(x)|\\
			\le&  |\widetilde{f}_{2,\Delta}(x)|\kappa_1^f(x_0)\Delta^{N_1^f} + |f_1(x)|\kappa_2^f(x_0)\Delta^{N_2^f}.
		\end{align*}
		Since $f_1$ and $f_2$ are interchangeable, we also have
		\begin{align*}
			|F(x) - \widetilde{F}_\Delta(x)| \le |f_2(x)|\kappa_1^f(x_0)\Delta^{N_1^f} + |\widetilde{f}_{1,\Delta}(x)|\kappa_2^f(x_0)\Delta^{N_2^f}.
		\end{align*}
		Therefore,
		\begin{align*}
			|F(x) - \widetilde{F}_\Delta(x)| \le& \min\{|\widetilde{f}_{2,\Delta}(x)|\kappa_1^f(x_0)\Delta^{N_1^f} + |f_1(x)|\kappa_2^f(x_0)\Delta^{N_2^f},\\
			&\qquad\ |f_2(x)|\kappa_1^f(x_0)\Delta^{N_1^f} + |\widetilde{f}_{1,\Delta}(x)|\kappa_2^f(x_0)\Delta^{N_2^f}\}\\
			\le& K_F^f(x_0)\Delta^{N_F^f}.
		\end{align*}
		That is, $\{\widetilde{F}_\Delta\}_{(0,\bar{\Delta}]}$ provides order-$N_F^f$ function accuracy of $F$ near $x_0$.
		
		$\hfill\qed$
	\end{proof}
	
	Theorem \ref{thm:2f} shows that, in terms of function accuracy, $\{\widetilde{F}_\Delta\}_{(0,\bar{\Delta}]}$ behaves at least as well as the worst behaved sub-model. Examining the result of Theorem \ref{thm:2f}, we can find a specific case where the overall function accuracy has a tighter bound. 
	
	\begin{corollary}\label{cor:2f}
		Suppose the conditions for Theorem \ref{thm:2f} hold. Without loss of generality, suppose $\min\{N_1^f,N_2^f\} = N_1^f$. If $f_2(x_0)=0$, then $\{\widetilde{F}_\Delta\}_{(0,\bar{\Delta}]}$ provides (at least) order-$N_2^f$ function accuracy of $F$ at $x_0$.
	\end{corollary}

	It should be noted that, if $f(x_0)$ is known, then given a class of models $\{\widetilde{f}_\Delta\}_{(0,\bar{\Delta}]}$ of $f$ we can construct a new class of models that obtains the true function value at $x_0$ without impacting gradient accuracy. Indeed, let
	\begin{align*}
		\widehat{f}_\Delta = \widetilde{f}_\Delta - \widetilde{f}_\Delta(x_0) + f(x_0).
	\end{align*}
	Then
	\begin{align*}
		\widehat{f}_\Delta(x_0) &= f(x_0),\\
		\nabla\widehat{f}_\Delta(x_0) &= \nabla\widetilde{f}_\Delta(x_0).
	\end{align*}

	Also notice that $\widehat{f}_\Delta(x_0) = f(x_0)$ for all $\Delta \in (0,\bar{\Delta}]$ is equivalent to saying that $\{\widehat{f}_\Delta\}_{(0,\bar{\Delta}]}$ provides order-$\infty$ function accuracy at $x_0$ with constant $\kappa^f(x_0)=0$. For the sake of consistency, we use order-$\infty$ terminology throughout this paper. Since we can apply the translation above to all the models, it is not unreasonable to assume that both $\{\widetilde{f}_{1,\Delta}\}_{(0,\bar{\Delta}]}$ and $\{\widetilde{f}_{2,\Delta}\}_{(0,\bar{\Delta}]}$ provide order-$\infty$ function accuracy at $x_0$ with constants $\kappa_1^f(x_0)=\kappa_2^f(x_0)=0$. The next theorem shows that this translation also creates a tighter bound for the overall gradient accuracy. Note, since we cannot suppose a class of models provides order-$\infty$ function accuracy near $x_0$ without loss of generality, we separate the ‘at $x_0$’ case and the ‘near $x_0$’ case for gradient and Hessian analysis.

	\begin{theorem}{(Gradient accuracy at $x_0$)}\label{thm:2ga}
		Let $f_1,f_2$ be in $\mathcal{C}^1$. Suppose for all $i\in\{1,2\}$, $\{\widetilde{f}_{i,\Delta}\}_{(0,\bar{\Delta}]}$ provides order-$\infty$ function accuracy of $f_i$ at $x_0$ with constant $\kappa_i^f(x_0)=0$, and order-$N_i^g$ gradient accuracy of $f_i$ at $x_0$ with constant $\kappa_i^g(x_0)$. Let $N_F^g=\min\{N_1^g,N_2^g\}$. Then
		\begin{align*}
			||\nabla F(x_0) - \nabla\widetilde{F}_\Delta(x_0)|| \le K_F^g(x_0)\Delta^{N_F^g},
		\end{align*}
		where $$K_F^g(x_0) = |f_2(x_0)|\kappa_1^g(x_0)\bar{\Delta}^{N_1^g - N_F^g} + |f_1(x_0)|\kappa_2^g(x_0)\bar{\Delta}^{N_2^g - N_F^g}.$$
		In particular, $\{\widetilde{F}_\Delta\}_{(0,\bar{\Delta}]}$ provides order-$N_F^g$ gradient accuracy of $F$ at $x_0$.
	\end{theorem}
	
	\begin{proof}
		Let $E_F^g = ||\nabla F(x_0)-\nabla\widetilde{F}_\Delta(x_0)||$. Then
		\begin{align*}
			E_F^g =& ||(f_2\nabla f_1 + f_1\nabla f_2 - \widetilde{f}_{2,\Delta}\nabla\widetilde{f}_{1,\Delta} - \widetilde{f}_{1,\Delta}\nabla\widetilde{f}_{2,\Delta})(x_0)||\\
			=& ||(f_2\nabla f_1-f_2\nabla\widetilde{f}_{1,\Delta}+f_1\nabla f_2-f_1\nabla\widetilde{f}_{2,\Delta} + f_2\nabla\widetilde{f}_{1,\Delta} -\widetilde{f}_{2,\Delta}\nabla\widetilde{f}_{1,\Delta}+f_1\nabla\widetilde{f}_{2,\Delta} -\widetilde{f}_{1,\Delta}\nabla\widetilde{f}_{2,\Delta})(x_0)||\\
			=& ||(f_1\nabla\widetilde{f}_{2,\Delta} -\widetilde{f}_{1,\Delta}\nabla\widetilde{f}_{2,\Delta} + f_2\nabla\widetilde{f}_{1,\Delta} -\widetilde{f}_{2,\Delta}\nabla\widetilde{f}_{1,\Delta} + f_2\nabla f_1-f_2\nabla\widetilde{f}_{1,\Delta}+f_1\nabla f_2-f_1\nabla\widetilde{f}_{2,\Delta})(x_0)||\\
			\le& ||\nabla\widetilde{f}_{2,\Delta}(x_0)||\kappa_1^f(x_0)\Delta^{N_1^f} + ||\nabla\widetilde{f}_{1,\Delta}(x_0)||\kappa_2^f(x_0)\Delta^{N_2^f} + |f_2(x_0)|\kappa_1^g(x_0)\Delta^{N_1^g} + |f_1(x_0)|\kappa_2^g(x_0)\Delta^{N_2^g}\\
			=& |f_2(x_0)|\kappa_1^g(x_0)\Delta^{N_1^g} + |f_1(x_0)|\kappa_2^g(x_0)\Delta^{N_2^g}\\
			\le& K_F^g(x_0)\Delta^{N_F^g}.
		\end{align*}
		That is, $\{\widetilde{F}_\Delta\}_{(0,\bar{\Delta}]}$ provides order-$N_F^g$ gradient accuracy of $F$ at $x_0$.
		
		$\hfill\qed$
	\end{proof}

    Theorem \ref{thm:2ga} shows that, in terms of gradient accuracy at $x_0$, $\{\widetilde{F}_\Delta\}_{(0,\bar{\Delta}]}$ behaves at least as well as the worst behaved sub-model. Also, there is a specific case where the overall gradient accuracy has a tighter bound. 
	
	\begin{corollary}\label{cor:2ga}
		Suppose the conditions for Theorem \ref{thm:2ga} hold. Without loss of generality, suppose $\min\{N_1^g,N_2^g\} = N_1^g$. If $f_2(x_0)=0$, then $\{\widetilde{F}_\Delta\}_{(0,\bar{\Delta}]}$ provides (at least) order-$N_2^g$ gradient accuracy of $F$ at $x_0$.\\
	\end{corollary}
	
	Since the order-$\infty$ function accuracy assumption does not hold without loss of generality for the ‘near $x_0$’ case, the result is weaker.
	
	\begin{theorem}{(Gradient accuracy near $x_0$)}\label{thm:2gn}
		Let $f_1,f_2$ be in $\mathcal{C}^1$. Suppose for all $i\in\{1,2\}$, $\{\widetilde{f}_{i,\Delta}\}_{(0,\bar{\Delta}]}$ provides order-$N_i^f$ function accuracy of $f_i$ near $x_0$, and order-$N_i^g$ gradient accuracy of $f_i$ near $x_0$. Then $\{\widetilde{F}_\Delta\}_{(0,\bar{\Delta}]}$ provides order-$N_F^{fg}$ gradient accuracy of $F$ near $x_0$, where $N_F^{fg} = \min\{N_1^f,N_2^f,N_1^g,N_2^g\}$.
	\end{theorem}
    	
	\begin{proof}
	    Let $E_F^g = ||\nabla F(x)-\nabla\widetilde{F}_\Delta(x)||$. Changing $x_0$ to $x$ in the proof of Theorem \ref{thm:2ga}, we obtain $$E_F^g \le ||\nabla\widetilde{f}_{2,\Delta}(x)||\kappa_1^f(x_0)\Delta^{N_1^f} + ||\nabla\widetilde{f}_{1,\Delta}(x)||\kappa_2^f(x_0)\Delta^{N_2^f} + |f_2(x)|\kappa_1^g(x_0)\Delta^{N_1^g} + |f_1(x)|\kappa_2^g(x_0)\Delta^{N_2^g}.$$
	    Since all the coefficients above are bounded, the proof is complete.
	    
	    $\hfill\qed$
	\end{proof}

    Theorem \ref{thm:2gn} shows that, in terms of gradient accuracy near $x_0$, $\{\widetilde{F}_\Delta\}_{(0,\bar{\Delta}]}$ behaves at least as well as the worst behaved sub-model, where both function and gradient behaviors are considered.
	
	The next example shows the necessity of $\{\widetilde{f}_{1,\Delta}\}_{(0,\bar{\Delta}]}$ and $\{\widetilde{f}_{2,\Delta}\}_{(0,\bar{\Delta}]}$ being uniformly bounded near $x_0$. In particular, if we only have the information of gradient accuracy, then the gradient error of $\{\widetilde{F}_\Delta\}_{(0,\bar{\Delta}]}$ can be infinitely large.

	\begin{example}
		Suppose for all $i\in\{1,2\}$, $\{\widetilde{f}_{i,\Delta}\}_{(0,\bar{\Delta}]}$ provides order-$N_i^g$ gradient accuracy of $f_i$ at $x_0$. If the uniformly bounded condition is not given, then $||\nabla F(x_0) - \nabla\widetilde{F}_\Delta(x_0)||$ can be infinitely large.
	\end{example}
	
	\begin{proof}
		Suppose $\nabla f_2(x_0)\neq 0$. Let $\widetilde{f}_{1,\Delta}(x)=f_1(x)+\frac{1}{\Delta}$ and $\widetilde{f}_{2,\Delta}(x)=f_2(x)$. Then, regardless of $N_1^g$ and $N_2^g$ we have
		\begin{align*}
			||\nabla f_1(x_0) - \nabla\widetilde{f}_{1,\Delta}(x_0)|| &= 0 \le \Delta^{N_1^g},\\
			||\nabla f_2(x_0) - \nabla\widetilde{f}_{2,\Delta}(x_0)|| &= 0 \le \Delta^{N_2^g},\\
			||\nabla F(x_0) - \nabla\widetilde{F}_\Delta(x_0)|| &= \frac{1}{\Delta}||\nabla f_2(x_0)||.
		\end{align*}
		Since $\lim\limits_{\Delta \to 0}\frac{1}{\Delta}||\nabla f_2(x_0)|| = \infty$, $||\nabla F(x_0) - \nabla\widetilde{F}_\Delta(x_0)||$ can be infinitely large.
		
		$\hfill\qed$
	\end{proof}

    With similar thoughts to function and gradient analysis, we also develop error bounds for Hessian. Again, we separate the ‘at $x_0$’ case and the ‘near $x_0$’ case.
    
	\begin{theorem}{(Hessian accuracy at $x_0$)}\label{thm:2ha}
		Let $f_1,f_2$ be in $\mathcal{C}^2$. Suppose for all $i\in\{1,2\}$,  $\{\widetilde{f}_{i,\Delta}\}_{(0,\bar{\Delta}]}$ provides order-$\infty$ function accuracy of $f_i$ at $x_0$ with constant $\kappa_i^f(x_0)=0$, order-$N_i^g$ gradient accuracy of $f_i$ at $x_0$ with constant $\kappa_i^g(x_0)$, and order-$N_i^h$ Hessian accuracy of $f_i$ at $x_0$ with constant $\kappa_i^h(x_0)$. Let $N_F^{gh}=\min\{N_1^g,N_2^g,N_1^h,N_2^h\}$. Then
		\begin{align*}
			||\nabla^2 F(x_0) - \nabla^2\widetilde{F}_\Delta(x_0)|| \le K_F^h(x_0)\Delta^{N_F^{gh}},
		\end{align*}
		where
		\begin{align*}
			K_F^h(x_0) =& (M_{\nabla f_2}(x_0) + M_{\nabla\widetilde{f}_2}(x_0))\kappa_1^g(x_0)\bar{\Delta}^{N_1^g - N_F^{gh}} + (M_{\nabla f_1}(x_0) + M_{\nabla\widetilde{f}_1}(x_0))\kappa_2^g(x_0)\bar{\Delta}^{N_2^g - N_F^{gh}}\\
			&+ |f_2(x_0)|\kappa_1^h(x_0)\bar{\Delta}^{N_1^h - N_F^{gh}} + |f_1(x_0)|\kappa_2^h(x_0)\bar{\Delta}^{N_2^h - N_F^{gh}}.
		\end{align*}
		In particular, $\{\widetilde{F}_\Delta\}_{(0,\bar{\Delta}]}$ provides order-$N_F^{gh}$ Hessian accuracy for $F$ at $x_0$.
	\end{theorem}
	
	\begin{proof}
		Let $E_F^h = ||\nabla^2 F(x_0) - \nabla^2\widetilde{F}_\Delta(x_0)||$. Then
		\begin{align*}
			E_F^h =& ||(f_2\nabla^2 f_1 + \nabla f_1\nabla f_2^\top + \nabla f_2\nabla f_1^\top + f_1\nabla^2 f_2 - \widetilde{f}_{2,\Delta}\nabla^2\widetilde{f}_{1,\Delta} - \nabla\widetilde{f}_{1,\Delta}\nabla\widetilde{f}_{2,\Delta}^\top - \nabla\widetilde{f}_{2,\Delta}\nabla\widetilde{f}_{1,\Delta}^\top\\
			&- \widetilde{f}_{1,\Delta}\nabla^2\widetilde{f}_{2,\Delta})(x_0)||\\
			=& ||(f_2\nabla^2 f_1 - f_2\nabla^2\widetilde{f}_{1,\Delta} + \nabla f_1\nabla f_2^\top - \nabla\widetilde{f}_{1,\Delta}\nabla f_2^\top + \nabla f_2\nabla f_1^\top - \nabla\widetilde{f}_{2,\Delta}\nabla f_1^\top + f_1\nabla^2 f_2\\ 
			&- f_1\nabla^2\widetilde{f}_{2,\Delta} - \widetilde{f}_{2,\Delta}\nabla^2\widetilde{f}_{1,\Delta} + f_2\nabla^2\widetilde{f}_{1,\Delta} - \nabla\widetilde{f}_{1,\Delta}\nabla\widetilde{f}_{2,\Delta}^\top + \nabla\widetilde{f}_{1,\Delta}\nabla f_2^\top - \nabla\widetilde{f}_{2,\Delta}\nabla\widetilde{f}_{1,\Delta}^\top\\
			&+ \nabla\widetilde{f}_{2,\Delta}\nabla f_1^\top - \widetilde{f}_{1,\Delta}\nabla^2\widetilde{f}_{2,\Delta} + f_1\nabla^2\widetilde{f}_{2,\Delta})(x_0)||\\
			=& ||(f_1\nabla^2\widetilde{f}_{2,\Delta} - \widetilde{f}_{1,\Delta}\nabla^2\widetilde{f}_{2,\Delta} + f_2\nabla^2\widetilde{f}_{1,\Delta} - \widetilde{f}_{2,\Delta}\nabla^2\widetilde{f}_{1,\Delta} + \nabla f_1\nabla f_2^\top - \nabla\widetilde{f}_{1,\Delta}\nabla f_2^\top + \nabla\widetilde{f}_{2,\Delta}\nabla f_1^\top\\
			&- \nabla\widetilde{f}_{2,\Delta}\nabla\widetilde{f}_{1,\Delta}^\top + \nabla f_2\nabla f_1^\top - \nabla\widetilde{f}_{2,\Delta}\nabla f_1^\top + \nabla\widetilde{f}_{1,\Delta}\nabla f_2^\top - \nabla\widetilde{f}_{1,\Delta}\nabla\widetilde{f}_{2,\Delta}^\top + f_2\nabla^2 f_1 - f_2\nabla^2\widetilde{f}_{1,\Delta}\\
			&+ f_1\nabla^2 f_2 - f_1\nabla^2\widetilde{f}_{2,\Delta})(x_0)||\\
			\le& ||\nabla^2\widetilde{f}_{2,\Delta}(x_0)||\kappa_1^f(x_0)\Delta^{N_1^f} + ||\nabla^2\widetilde{f}_{1,\Delta}(x_0)||\kappa_2^f(x_0)\Delta^{N_2^f} + (||\nabla f_2(x_0)|| + ||\nabla\widetilde{f}_{2,\Delta}(x_0)||) \kappa_1^g(x_0)\Delta^{N_1^g}\\
			&+ (||\nabla f_1(x_0)|| + ||\nabla\widetilde{f}_{1,\Delta}(x_0)||)\kappa_2^g(x_0)\Delta^{N_2^g} + |f_2(x_0)|\kappa_1^h(x_0)\Delta^{N_1^h} + |f_1(x_0)|\kappa_2^h(x_0)\Delta^{N_2^h}\\
			=& (||\nabla f_2(x_0)|| + ||\nabla\widetilde{f}_{2,\Delta}(x_0)||) \kappa_1^g(x_0)\Delta^{N_1^g} + (||\nabla f_1(x_0)|| + ||\nabla\widetilde{f}_{1,\Delta}(x_0)||)\kappa_2^g(x_0)\Delta^{N_2^g}\\
			&+ |f_2(x_0)|\kappa_1^h(x_0)\Delta^{N_1^h} + |f_1(x_0)|\kappa_2^h(x_0)\Delta^{N_2^h}\\
			\le& K_F^h(x_0)\Delta^{N_F^{gh}}.
		\end{align*}
		
		$\hfill\qed$
	\end{proof}

    Theorem \ref{thm:2ha} shows that, in terms of Hessian accuracy at $x_0$, $\{\widetilde{F}_\Delta\}_{(0,\bar{\Delta}]}$ behaves at least as well as the worst behaved sub-model, where both gradient and Hessian behaviors are considered. 
	
	\begin{theorem}{(Hessian accuracy near $x_0$)}\label{thm:2hn}
		Let $f_1,f_2$ be in $\mathcal{C}^2$. Suppose for all $i\in\{1,2\}$, $\{\widetilde{f}_{i,\Delta}\}_{(0,\bar{\Delta}]}$ provides order-$N_i^f$ function accuracy of $f_i$ near $x_0$, order-$N_i^g$ gradient accuracy of $f_i$ near $x_0$, and order-$N_i^h$ Hessian accuracy of $f_i$ near $x_0$. Then $\{\widetilde{F}_\Delta\}_{(0,\bar{\Delta}]}$ provides order-$N_F^{fgh}$ Hessian accuracy for $F$ near $x_0$, where $N_F^{fgh}=\min\{N_1^f,N_2^f,N_1^g,N_2^g,N_1^h,N_2^h\}$.
	\end{theorem}
	
	\begin{proof}
	    Let $E_F^h = ||\nabla^2 F(x) - \nabla^2\widetilde{F}_\Delta(x)||$. Changing $x_0$ to $x$ in the proof of Theorem \ref{thm:2ha}, we obtain
	    \begin{align*}
	        E_F^h \le& ||\nabla^2\widetilde{f}_{2,\Delta}(x)||\kappa_1^f(x_0)\Delta^{N_1^f} + ||\nabla^2\widetilde{f}_{1,\Delta}(x)||\kappa_2^f(x_0)\Delta^{N_2^f} + (||\nabla f_2(x)|| + ||\nabla\widetilde{f}_{2,\Delta}(x)||) \kappa_1^g(x_0)\Delta^{N_1^g}\\
			&+ (||\nabla f_1(x)|| + ||\nabla\widetilde{f}_{1,\Delta}(x)||)\kappa_2^g(x_0)\Delta^{N_2^g} + |f_2(x)|\kappa_1^h(x_0)\Delta^{N_1^h} + |f_1(x)|\kappa_2^h(x_0)\Delta^{N_2^h}.
	    \end{align*}
	    Since all the coefficients above are bounded, the proof is complete.
	    
	    $\hfill\qed$
	\end{proof}

    Theorem \ref{thm:2hn} shows that, in terms of Hessian accuracy near $x_0$, $\{\widetilde{F}_\Delta\}_{(0,\bar{\Delta}]}$ behaves at least as well as the worst behaved sub-model, where all function, gradient, and Hessian behaviors are considered.

	In the next subsection, we extend the results in this subsection to $n$-function case.
	\subsection{Error Analysis of the Product of $n$ Functions}
	
	Throughout this subsection, we use $F=f_1\cdots f_n$ and $\widetilde{F}_\Delta=\widetilde{f}_{1,\Delta}\cdots\widetilde{f}_{n,\Delta}$. In the next theorem, we use $S_n$ as the permutation group of the set $\{1,...,n\}$, and the function $\sigma$ as any element of $S_n$.
	
	\begin{theorem}{(Function accuracy at/near $x_0$)}\label{thm:nf}
		For all $i\in\{1,...,n\}$, let $f_i$ be in $\mathcal{C}^0$. Suppose $\{\widetilde{f}_{i,\Delta}\}_{(0,\bar{\Delta}]}$ provides order-$N_i^f$ function accuracy of $f_i$ at/near $x_0$ with constant $\kappa_i^f(x_0)$. Let \\$N_F^f = \min\{N_i^f:i=1,...,n\}$. Then
		
		$\ $(i) if the accuracy provided is at $x_0$, then
		\begin{align*}
			|F(x_0) - \widetilde{F}_\Delta(x_0)| \le K_F^f(x_0)\Delta^{N_F^f};
		\end{align*}
	
		(ii) if the accuracy provided is near $x_0$, then for all $x\in B_\Delta(x_0)$
		\begin{align*}
			|F(x) - \widetilde{F}_\Delta(x)| \le K_F^f(x_0)\Delta^{N_F^f};
		\end{align*}
		where 
		\begin{align}\label{eq:K^f}
			K_F^f(x_0) = \min\left\{\sum\limits_{i=1}^n\left(\prod\limits_{j=1}^{i-1}M_{f_{\sigma(j)}}(x_0)\prod\limits_{k=i+1}^n M_{\widetilde{f}_{\sigma(k)}}(x_0)\right)\kappa_{\sigma(i)}^f(x_0)\bar{\Delta}^{N_{\sigma(i)}^f-N_F^f}:\sigma\in S_n\right\}.
		\end{align}
		In particular, $\{\widetilde{F}_\Delta\}_{(0,\bar{\Delta}]}$ provides order-$N_F^f$ function accuracy of $F$ at/near $x_0$.
	\end{theorem}
	
	\begin{proof}
		We only prove case (ii), since the proof of case (i) can easily be obtained by changing all $x$ to $x_0$. By induction from Theorem \ref{thm:2f}, we have
		\begin{align*}
			|F(x)-\widetilde{F}_\Delta(x)|\le\sum\limits_{i=1}^n\left(\prod\limits_{j=1}^{i-1}|f_j(x)|\prod\limits_{k=i+1}^n|\widetilde{f}_{k,\Delta}(x)|\right)\kappa_i^f(x_0)\Delta^{N_i^f}.
		\end{align*}
		Since each two $f_i$ are interchangeable, we also have
		\begin{align*}
			|F(x)-\widetilde{F}_\Delta(x)|\le\sum\limits_{i=1}^n\left(\prod\limits_{j=1}^{i-1}|f_{\sigma(j)}(x)|\prod\limits_{k=i+1}^n|\widetilde{f}_{\sigma(k),\Delta}(x)|\right)\kappa_{\sigma(i)}^f(x_0)\Delta^{N_{\sigma(i)}^f}.
		\end{align*}
		Therefore,
		\begin{align*}
			|F(x) - \widetilde{F}_\Delta(x)| \le& \min\left\{\sum\limits_{i=1}^n\left(\prod\limits_{j=1}^{i-1}|f_{\sigma(j)}(x)|\prod\limits_{k=i+1}^n|\widetilde{f}_{\sigma(k),\Delta}(x)|\right)\kappa_{\sigma(i)}^f(x_0)\Delta^{N_{\sigma(i)}^f}:\sigma\in S_n\right\}\\
			\le& K_F^f(x_0)\Delta^{N_F^f}.
		\end{align*}
		That is, $\{\widetilde{F}_\Delta\}_{(0,\bar{\Delta}]}$ provides order-$N_F^f$ function accuracy of $F$ near $x_0$.
		
		$\hfill\qed$
	\end{proof}

	Theorem \ref{thm:nf} shows that, in terms of function accuracy, $\{\widetilde{F}_\Delta\}_{(0,\bar{\Delta}]}$ behaves at least as well as the worst behaved sub-model. Examining the result of Theorem \ref{thm:nf}, we can find a specific case where the overall function accuracy has a tighter bound. 
	
	\begin{corollary}
		Suppose the conditions for Theorem \ref{thm:nf} hold. If there exist $i\in\{1,...,n\}$ such that $f_i(x_0)=0$, then $\{\widetilde{F}_\Delta\}_{(0,\bar{\Delta}]}$ provides (at least) order-$N_i^f$ function accuracy of $F$ at $x_0$.
	\end{corollary}

	With similar thought to function accuracy analysis, we also develop error bounds for gradient and Hessian.
	
	\begin{theorem}{(Gradient accuracy at $x_0$)}\label{thm:nga}
		For all $i\in\{1,...,n\}$, let $f_i$ be in $\mathcal{C}^1$. Suppose $\{\widetilde{f}_{i,\Delta}\}_{(0,\bar{\Delta}]}$ provides order-$\infty$ function accuracy of $f_i$ at $x_0$ with constant $\kappa_i^f(x_0)=0$, and order-$N_i^g$ gradient accuracy of $f_i$ at $x_0$ with constant $\kappa_i^g(x_0)$. Let $N_F^g = \min\{N_i^g:i=1,...,n\}$. Then
		\begin{align*}
			||\nabla F(x_0) - \nabla\widetilde{F}_\Delta(x_0)|| \le K_F^g(x_0)\Delta^{N_F^g},
		\end{align*}
		where 
		\begin{align}\label{eq:K^g}
		    K_F^g = \sum\limits_{i=1}^n\left(\prod\limits_{j\neq i}^n|f_j(x_0)|\right)\kappa_i^g(x_0)\bar{\Delta}^{N_i^g - N_F^g}.
		\end{align} 
		In particular, $\{\widetilde{F}_\Delta\}_{(0,\bar{\Delta}]}$ provides order-$N_F^g$ gradient accuracy of $F$ at $x_0$.
	\end{theorem}
	
	\begin{proof}
		Let $F_i=f_1\cdots f_{i-1}f_{i+1}\cdots f_n,\widetilde{F}_{i,\Delta}=\widetilde{f}_{1,\Delta}\cdots \widetilde{f}_{i-1,\Delta}\widetilde{f}_{i+1,\Delta}\cdots \widetilde{f}_{n,\Delta}$. Let $E_F^g = ||\nabla F(x_0) - \nabla\widetilde{F}_\Delta(x_0)||$. By induction from Theorem \ref{thm:2ga}, we have
		\begin{align*}
			E_F^g \le&
			\sum\limits_{i=1}^n\left[||\nabla\widetilde{f}_{i,\Delta}(x_0)|||F_i(x_0)-\widetilde{F}_{i,\Delta}(x_0)| + \left(\prod\limits_{j\neq i}^n|f_j(x_0)|\right)\kappa_i^g(x_0)\Delta^{N_i^g}\right]\\
			=& \sum\limits_{i=1}^n\left(\prod\limits_{j\neq i}^n|f_j(x_0)|\right)\kappa_i^g(x_0)\Delta^{N_i^g}\\
			\le& K_F^g(x_0)\Delta^{N_F^g}.
		\end{align*}
		That is, $\{\widetilde{F}_\Delta\}_{(0,\bar{\Delta}]}$ provides order-$N_F^g$ gradient accuracy of $F$ at $x_0$.
		
		$\hfill\qed$
	\end{proof}
	
	Theorem \ref{thm:nga} shows that, in terms of gradient accuracy at $x_0$, $\{\widetilde{F}_\Delta\}_{(0,\bar{\Delta}]}$ behaves at least as well as the worst behaved sub-model. Examining the result of Theorem \ref{thm:nga}, we find cases where the overall gradient accuracy has a tighter bound. It is worth noticing that if any two different component functions are zero at $x_0$, then we will get the true gradient of $F$ at $x_0$.
	
	\begin{corollary}\label{cor:nga}
		Suppose the conditions for Theorem \ref{thm:nga} hold. If there exists $i\in\{1,...,n\}$ such that $f_i(x_0)=0$, then $\{\widetilde{F}_\Delta\}_{(0,\bar{\Delta}]}$ provides (at least) order-$N_i^g$ gradient accuracy of $F$ at $x_0$. If there exists $i,j\in\{1,...,n\},i\neq j$ such that $f_i(x_0)=f_j(x_0)=0$, then $\{\widetilde{F}_\Delta\}_{(0,\bar{\Delta}]}$ obtains the true gradient of $F$ at $x_0$.
	\end{corollary}

    Since the order-$\infty$ function accuracy assumption does not hold without loss of generality for the ‘near $x_0$’ case, the result is weaker.
	
	\begin{theorem}{(Gradient accuracy near $x_0$)}\label{thm:ngn}
		For all $i\in\{1,...,n\}$, let $f_i$ be in $\mathcal{C}^1$. Suppose $\{\widetilde{f}_{i,\Delta}\}_{(0,\bar{\Delta}]}$ provides order-$N_i^f$ function accuracy of $f_i$ near $x_0$, and order-$N_i^g$ gradient accuracy of $f_i$ near $x_0$. Then $\{\widetilde{F}_\Delta\}_{(0,\bar{\Delta}]}$ provides order-$N_F^{fg}$ gradient accuracy of $F$ near $x_0$, where $N_F^{fg} = \min\{N_i^f,N_i^g:i=1,...,n\}$.
	\end{theorem}
	
	\begin{proof}
	    Let $E_F^g = ||\nabla F(x)-\nabla\widetilde{F}_\Delta(x)||$. Changing $x_0$ to $x$ in the proof of Theorem \ref{thm:nga}, we obtain $$E_F^g \le
		\sum\limits_{i=1}^n\left[||\nabla\widetilde{f}_{i,\Delta}(x)|||F_i(x)-\widetilde{F}_{i,\Delta}(x)| + \left(\prod\limits_{j\neq i}^n|f_j(x)|\right)\kappa_i^g(x_0)\Delta^{N_i^g}\right].$$
		Notice that, since the term $|F_i(x)-\widetilde{F}_{i,\Delta}(x)|$ can be bounded by some weighted sum of $\kappa_i^f(x_0)\Delta^{N_i^f}$ with all the coefficients bounded, the proof is complete.
		
		$\hfill\qed$
	\end{proof}

    Theorem \ref{thm:ngn} shows that, in terms of gradient accuracy near $x_0$, $\{\widetilde{F}_\Delta\}_{(0,\bar{\Delta}]}$ behaves at least as well as the worst behaved sub-model, where both function and gradient behaviors are considered. 

	\begin{theorem}{(Hessian accuracy at $x_0$)}\label{thm:nha}
		For all $i\in\{1,...,n\}$, let $f_i$ be in $\mathcal{C}^2$. Suppose $\{\widetilde{f}_{i,\Delta}\}_{(0,\bar{\Delta}]}$ provides order-$\infty$ function accuracy of $f_i$ at $x_0$ with constant $\kappa_i^f(x_0)=0$, order-$N_i^g$ gradient accuracy of $f_i$ at $x_0$ with constant $\kappa_i^g(x_0)$, and order-$N_i^h$ Hessian accuracy of $f_i$ at $x_0$ with constant $\kappa_i^h(x_0)$. Let $N_F^{gh}=\min\{N_i^g,N_i^h:i=1,...,n\}$. Then
		\begin{align*}
			||\nabla^2 F(x_0) - \nabla^2\widetilde{F}_\Delta(x_0)|| \le K_F^h(x_0)\Delta^{N_F^{gh}},
		\end{align*}
		where
		\begin{equation}\label{eq:K^h}
		    \begin{aligned}
    			K_F^h(x_0) =& \sum\limits_{i=1}^n\left[ \sum\limits_{j\neq i}^n\left( \prod\limits_{k\neq i,j}^n|f_k(x_0)| \right)(M_{\nabla f_j}(x_0) + M_{\nabla\widetilde{f}_j}(x_0))\right]\kappa_i^g(x_0)\bar{\Delta}^{N_i^g - N_F^{gh}}\\
    			&+ \sum\limits_{i=1}^n |F_i(x_0)|\kappa_i^h(x_0)\bar{\Delta}^{N_i^h - N_F^{gh}}.
    		\end{aligned}
		\end{equation}
		In particular, $\{\widetilde{F}_\Delta\}_{(0,\bar{\Delta}]}$ provides order-$N_F^{gh}$ Hessian accuracy for $F$ at $x_0$.
	\end{theorem}
	
    \begin{proof}
		Let $E_F^h = ||\nabla^2 F(x_0) - \nabla^2\widetilde{F}_\Delta(x_0)||$. By induction from Theorem \ref{thm:2ha}, we have
		\begin{align*}
			E_F^h \le& \sum\limits_{i=1}^n\left[ ||\nabla^2\widetilde{f}_{i,\Delta}(x_0)|||F_i(x_0) - \widetilde{F}_{i,\Delta}(x_0)| + ||\nabla\widetilde{f}_{i,\Delta}(x_0)||||\nabla F_i(x_0) - \nabla\widetilde{F}_{i,\Delta}(x_0)|| \right.\\
			&\left.+ ||\nabla F_i(x_0)||\kappa_i^g(x_0)\Delta^{N_i^g} + |F_i(x_0)|\kappa_i^h(x_0)\Delta^{N_i^h} \right]\\
			\le& \sum\limits_{i=1}^n\left\{ \left[||\nabla F_i(x_0)|| + \sum\limits_{j\neq i}^n\left(\prod\limits_{k\neq i,j}^n|f_k(x_0)|\right)||\nabla\widetilde{f}_{i,\Delta}(x_0)||\right] \kappa_i^g(x_0)\Delta^{N_i^g} + |F_i(x_0)|\kappa_i^h(x_0)\Delta^{N_i^h} \right\}\\
			=& \sum\limits_{i=1}^n\left[ \sum\limits_{j\neq i}^n\left( \prod\limits_{k\neq i,j}^n|f_k(x_0)| \right)(||\nabla f_j(x_0)|| + ||\nabla\widetilde{f}_j(x_0)||)\right]\kappa_i^g(x_0)\Delta^{N_i^g} + \sum\limits_{i=1}^n |F_i(x_0)|\kappa_i^h(x_0)\Delta^{N_i^h}\\
			\le& K_F^h(x_0)\Delta^{N_F^{gh}}.
		\end{align*}
		
		$\hfill\qed$
	\end{proof}

    Theorem \ref{thm:nha} shows that, in terms of Hessian accuracy at $x_0$, $\{\widetilde{F}_\Delta\}_{(0,\bar{\Delta}]}$ behaves at least as well as the worst behaved sub-model, where both gradient and Hessian behaviors are considered.

	\begin{theorem}{(Hessian accuracy near $x_0$)}\label{thm:nhn}
		For all $i\in\{1,...,n\}$, let $f_i$ be in $\mathcal{C}^2$. Suppose $\{\widetilde{f}_{i,\Delta}\}_{(0,\bar{\Delta}]}$ provides order-$N_i^f$ function accuracy of $f_i$ near $x_0$, order-$N_i^g$ gradient accuracy of $f_i$ near $x_0$, and order-$N_i^h$ Hessian accuracy of $f_i$ near $x_0$. Then $\{\widetilde{F}_\Delta\}_{(0,\bar{\Delta}]}$ provides order-$N_F^{fgh}$ Hessian accuracy for $F$ near $x_0$, where $N_F^{fgh}=\min\{N_i^f,N_i^g,N_i^h:i=1,...,n\}$.
	\end{theorem}
	
	\begin{proof}
	    Let $E_F^h = ||\nabla^2 F(x)-\nabla^2\widetilde{F}_\Delta(x)||$. Changing $x_0$ to $x$ in the proof of Theorem \ref{thm:nha}, we obtain $$E_F^h\le\sum\limits_{i=1}^n\left\{ \left[||\nabla F_i(x)|| + \sum\limits_{j\neq i}^n\left(\prod\limits_{k\neq i,j}^n|f_k(x)|\right)||\nabla\widetilde{f}_{i,\Delta}(x)||\right] \kappa_i^g(x_0)\Delta^{N_i^g} + |F_i(x)|\kappa_i^h(x_0)\Delta^{N_i^h} \right\}.$$
	    Since all the coefficients above are bounded, the proof is complete.
	    
	    $\hfill\qed$
	\end{proof}
	
	Theorem \ref{thm:nhn} shows that, in terms of Hessian accuracy near $x_0$, $\{\widetilde{F}_\Delta\}_{(0,\bar{\Delta}]}$ behaves at least as well as the worst behaved sub-model, where all function, gradient, and Hessian behaviors are considered.	
	
	Note that, simply letting $f_1=...=f_n=f$ in Theorem \ref{thm:nf}, Theorem \ref{thm:nga}, Theorem \ref{thm:ngn}, Theorem \ref{thm:nha}, and Theorem \ref{thm:nhn}, we obtain error bounds for the case where the objective function can be expressed as the power of a component function. Specifically, we obtain $$K_F^f(x_0) = \left(\sum\limits_{i=1}^n M_f(x_0)^{i-1}M_{\widetilde{f}}(x_0)^{n-i}\right)\kappa^f(x_0),$$ $$K_F^g(x_0) = n|f(x_0)|^{n-1}\kappa^g(x_0),$$ and $$K_F^h(x_0) = n(n-1)|f(x_0)|^{n-2}(M_{\nabla f}(x_0) + M_{\nabla\widetilde{f}}(x_0))\kappa^g(x_0)\bar{\Delta}^{N^g - N_F^{gh}} + n|f(x_0)|^{n-1}\kappa^h(x_0)\bar{\Delta}^{N^h - N_F^{gh}}$$ from (\ref{eq:K^f}), (\ref{eq:K^g}), and (\ref{eq:K^h}), respectively.

	\section{Error Analysis of the Quotient of Model Functions}
	
	In this section, we consider the case where the objective function $F$ can be expressed as the quotient of 2 component functions. Throughout this section, suppose $f_2(x)\neq 0,\widetilde{f}_{2,\Delta}(x)\neq 0$ for all $\Delta\in(0,\bar{\Delta}]$ and $x\in B_\Delta(x_0)$. Let $F=\frac{f_1}{f_2},\widetilde{F}_\Delta=\frac{\widetilde{f}_{1,\Delta}}{\widetilde{f}_{2,\Delta}}$.
	
	From fundamental mathematical analysis we have that if the function $f_2\in\mathcal{C}^0$ and $f_2(x_0)\neq 0$, then we can choose a sufficiently small $\bar{\Delta}$ such that $\frac{1}{f_2}$ is bounded in $B_{\bar{\Delta}}(x_0)$. We denote this bound by $M_{\frac{1}{f_2}}(x_0)$. The next lemma shows that under certain conditions, $\{\frac{1}{\widetilde{f}_{2,\Delta}}\}_{(0,\bar{\Delta}]}$ is uniformly bounded.
	
	\begin{lemma}\label{lem:ub}
		Let $f_2$ be in $\mathcal{C}^0$. Suppose $\{\widetilde{f}_{2,\Delta}\}_{(0,\bar{\Delta}]}$ provides order-$N_2^f$ function accuracy of $f_2$ at/near $x_0$ with constant $\kappa_2^f(x_0)$. If $N_2^f > 0$, then we can choose a sufficiently small $\bar{\Delta}$ such that $\{\frac{1}{\widetilde{f}_{2,\Delta}}\}_{(0,\bar{\Delta}]}$ is uniformly bounded at/near $x_0$ by $$\widetilde{M}_{\frac{1}{\widetilde{f}_2}}(x_0) = \frac{M_{\frac{1}{f_2}}(x_0)}{1-\bar{\Delta}^{N_2^f}M_{\frac{1}{f_2}}(x_0)\kappa_2^f(x_0)}.$$
	\end{lemma}

	\begin{proof}
		We only prove the ‘near $x_0$’ case, since the ‘at $x_0$’ case is similar.\\
		For all $\Delta\in(0,\bar{\Delta}]$ and $x\in B_{\Delta}(x_0)$,
		\begin{align*}
			\left|\frac{1}{\widetilde{f}_{2,\Delta}(x)}\right| \le& \left|\frac{1}{\widetilde{f}_{2,\Delta}(x)} - \frac{1}{f_2(x)}\right| + \left|\frac{1}{f_2(x)}\right|\\
			\le& \left|\frac{f_2(x) - \widetilde{f}_{2,\Delta}(x)}{f_2(x)\widetilde{f}_{2,\Delta}(x)}\right| + M_{\frac{1}{f_2}}(x_0)\\
			\le& M_{\frac{1}{f_2}}(x_0)\kappa_2^f(x_0)\Delta^{N_2^f}\left|\frac{1}{\widetilde{f}_{2,\Delta}(x)}\right| + M_{\frac{1}{f_2}}(x_0)\\
			\le& \bar{\Delta}^{N_2^f}M_{\frac{1}{f_2}}(x_0)\kappa_2^f(x_0)\left|\frac{1}{\widetilde{f}_{2,\Delta}(x)}\right| + M_{\frac{1}{f_2}}(x_0).
		\end{align*}
		If $N_2^f > 0$, we can choose a sufficiently small $\bar{\Delta}$ such that $1-\bar{\Delta}^{N_2^f}M_{\frac{1}{f_2}}(x_0)\kappa_2^f(x_0) > 0$. Therefore, we get
		\begin{align*}
			\left|\frac{1}{\widetilde{f}_{2,\Delta}(x)}\right| \le \widetilde{M}_{\frac{1}{\widetilde{f}_2}}(x_0).
		\end{align*}
		
		$\hfill\qed$
	\end{proof}

	Using the two bounds above, we develop the overall function error bound.
	
	\begin{theorem}{(Function accuracy at/near $x_0$)}\label{thm:qf}
		Let $f_1,f_2$ be in $\mathcal{C}^0$. Suppose for all $i\in\{1,2\}$, $\{\widetilde{f}_{i,\Delta}\}_{(0,\bar{\Delta}]}$ provides order-$N_i^f$ function accuracy of $f_i$ at/near $x_0$ with constant $\kappa_i^f(x_0)$. Let $N_F^f = \min\{N_1^f,N_2^f\}$. If $N_2^f>0$. Then
		
		$\ $(i) if the accuracy provided is at $x_0$, then
		\begin{align*}
			|F(x_0) - \widetilde{F}_\Delta(x_0)| \le K_F^f(x_0)\Delta^{N_F^f};
		\end{align*}
	
		(ii) if the accuracy provided is near $x_0$, then for all $x\in B_\Delta(x_0)$
		\begin{align*}
			|F(x) - \widetilde{F}_\Delta(x)| \le K_F^f(x_0)\Delta^{N_F^f};
		\end{align*}
		where 
		\begin{align*}
			K_F^f(x_0) = \min&\left\{M_{\frac{1}{f_2}}(x_0)\kappa_1^f(x_0)\bar{\Delta}^{N_1^f-N_F^f} + 
			M_{\frac{1}{f_2}}(x_0)\widetilde{M}_{\frac{1}{\widetilde{f}_2}}(x_0)M_{\widetilde{f}_1}(x_0)\kappa_2^f(x_0)\bar{\Delta}^{N_2^f-N_F^f},\right.\\
			&\ \ \left.\widetilde{M}_{\frac{1}{\widetilde{f}_2}}(x_0)\kappa_1^f(x_0)\bar{\Delta}^{N_1^f-N_F^f}   + M_{\frac{1}{f_2}}(x_0)\widetilde{M}_{\frac{1}{\widetilde{f}_2}}(x_0)M_{f_1}(x_0)\kappa_2^f(x_0)\bar{\Delta}^{N_2^f-N_F^f} \right\}.
		\end{align*}
		In particular, $\{\widetilde{F}_\Delta\}_{(0,\bar{\Delta}]}$ provides order-$N_F^f$ function accuracy of $F$ at/near $x_0$.
	\end{theorem}

	\begin{proof}
		We only prove case (ii), since the proof of case (i) can easily be obtained by changing all $x$ to $x_0$. First notice
		\begin{align*}
			|F(x) - \widetilde{F}_\Delta(x)| =& \left|\left(\frac{f_1}{f_2} - \frac{\widetilde{f}_{1,\Delta}}{\widetilde{f}_{2,\Delta}}\right)(x)\right|\\
			=& \left|\left(\frac{f_1}{f_2} - \frac{\widetilde{f}_{1,\Delta}}{f_2} + \frac{\widetilde{f}_{1,\Delta}}{f_2} - \frac{\widetilde{f}_{1,\Delta}}{\widetilde{f}_{2,\Delta}}\right)(x)\right|\\
			\le& \left|\frac{1}{f_2(x)}\right|\kappa_1^f(x_0)\Delta^{N_1^f} + \left|\frac{\widetilde{f}_{1,\Delta}(x)}{f_2(x)\widetilde{f}_{2,\Delta}(x)}\right|\kappa_2^f(x_0)\Delta^{N_2^f}.
		\end{align*}
		Similarly, we also have
		\begin{align*}
			|F(x) - \widetilde{F}_\Delta(x)| \le& \left|\frac{1}{\widetilde{f}_{2,\Delta}(x)}\right|\kappa_1^f(x_0)\Delta^{N_1^f} + \left|\frac{f_1(x)}{f_2(x)\widetilde{f}_{2,\Delta}(x)}\right|\kappa_2^f(x_0)\Delta^{N_2^f}.
		\end{align*}
		Therefore,
		\begin{align*}
			|F(x) - \widetilde{F}_\Delta(x)| \le&  \min\left\{\left|\frac{1}{f_2(x)}\right|\kappa_1^f(x_0)\Delta^{N_1^f} + 
			\left|\frac{\widetilde{f}_{1,\Delta}(x)}{f_2(x)\widetilde{f}_{2,\Delta}(x)}\right|\kappa_2^f(x_0)\Delta^{N_2^f},\right.\\
			&\qquad\  \left.\left|\frac{1}{\widetilde{f}_{2,\Delta}(x)}\right|\kappa_1^f(x_0)\Delta^{N_1^f} + \left|\frac{f_1(x)}{f_2(x)\widetilde{f}_{2,\Delta}(x)}\right|\kappa_2^f(x_0)\Delta^{N_2^f} \right\}\\
			\le& K_F^f(x_0)\Delta^{N_F^f}.
		\end{align*}
		That is, $\{\widetilde{F}_\Delta\}_{(0,\bar{\Delta}]}$ provides order-$N_F^f$ function accuracy of $F$ near $x_0$.
		
		$\hfill\qed$
	\end{proof}
	
	Theorem \ref{thm:qf} shows that, in terms of function accuracy, $\{\widetilde{F}_\Delta\}_{(0,\bar{\Delta}]}$ behaves at least as well as the worst behaved sub-model.
	
	Similarly, we also develop error bounds for gradient and Hessian.
	
	\begin{theorem}{(Gradient accuracy at $x_0$)}\label{thm:qga}
		Let $f_1,f_2$ be in $\mathcal{C}^1$. Suppose for all $i\in\{1,2\}$, $\{\widetilde{f}_{i,\Delta}\}_{(0,\bar{\Delta}]}$ provides order-$\infty$ function accuracy of $f_i$ at $x_0$ with constant $\kappa_i^f(x_0)=0$, and order-$N_i^g$ gradient accuracy of $f_i$ at $x_0$ with constant $\kappa_i^g(x_0)$. Let $N_F^g = \min\{N_1^g,N_2^g\}$. Then
		\begin{align*}
			||\nabla F(x_0) - \nabla\widetilde{F}_\Delta(x_0)|| \le K_F^g(x_0)\Delta^{N_F^g},
		\end{align*}
		where $$K_F^g(x_0) = \left|\frac{1}{f_2(x_0)}\right|\kappa_1^g(x_0)\bar{\Delta}^{N_1^g - N_F^g} + \left|\frac{f_1(x_0)}{f_2^2(x_0)}\right|\kappa_2^g(x_0)\bar{\Delta}^{N_2^g - N_F^g}.$$
		In particular, $\{\widetilde{F}_\Delta\}_{(0,\bar{\Delta}]}$ provides order-$N_F^g$ gradient accuracy of $F$ at $x_0$.
	\end{theorem}

	\begin{proof}
		Let $E_F^g = ||\nabla F(x_0) - \nabla\widetilde{F}_\Delta(x_0)||$. Then
		\begin{align*}
			E_F^g =& \left\|\left(\frac{\nabla f_1}{f_2} - \frac{f_1}{f_2^2}\nabla f_2 - \frac{\nabla\widetilde{f}_{1,\Delta}}{\widetilde{f}_{2,\Delta}} + \frac{\widetilde{f}_{1,\Delta}}{\widetilde{f}_{2,\Delta}^2}\nabla\widetilde{f}_{2,\Delta}\right)(x_0)\right\|\\
			=& \left\|\left(\frac{\nabla f_1}{f_2} - \frac{\nabla\widetilde{f}_{1,\Delta}}{f_2} - \frac{f_1}{f_2^2}\nabla f_2 + \frac{f_1}{f_2^2}\nabla\widetilde{f}_{2,\Delta} + \frac{\nabla\widetilde{f}_{1,\Delta}}{f_2} - \frac{\nabla\widetilde{f}_{1,\Delta}}{\widetilde{f}_{2,\Delta}} - \frac{f_1}{f_2^2}\nabla\widetilde{f}_{2,\Delta} + \frac{\widetilde{f}_{1,\Delta}}{\widetilde{f}_{2,\Delta}^2}\nabla\widetilde{f}_{2,\Delta}\right)(x_0)\right\|\\
			=& \left\|\left(\frac{\widetilde{f}_{1,\Delta}}{\widetilde{f}_{2,\Delta}^2}\nabla\widetilde{f}_{2,\Delta} - \frac{f_1}{f_2^2}\nabla\widetilde{f}_{2,\Delta} + \frac{\nabla\widetilde{f}_{1,\Delta}}{f_2} - \frac{\nabla\widetilde{f}_{1,\Delta}}{\widetilde{f}_{2,\Delta}} + \frac{\nabla f_1}{f_2} - \frac{\nabla\widetilde{f}_{1,\Delta}}{f_2} + \frac{f_1}{f_2^2}\nabla\widetilde{f}_{2,\Delta} - \frac{f_1}{f_2^2}\nabla f_2\right)(x_0)\right\|\\
			\le& \left\|\frac{\nabla\widetilde{f}_{2,\Delta}(x_0)}{f_2^2(x_0)\widetilde{f}_{2,\Delta}^2(x_0)}\right\|\left[\left|f_1(x_0)\left(f_2(x_0) + \widetilde{f}_{2,\Delta}(x_0)\right)\right|\kappa_2^f(x_0)\Delta^{N_2^f} + |f_2^2(x_0)|\kappa_1^f(x_0)\Delta^{N_1^f}\right]\\
			&+ \left\|\frac{\nabla \widetilde{f}_{1,\Delta}(x_0)}{f_2(x_0)\widetilde{f}_{2,\Delta}(x_0)}\right\|\kappa_2^f(x_0)\Delta^{N_2^f} + \left|\frac{1}{f_2(x_0)}\right|\kappa_1^g(x_0)\Delta^{N_1^g} + \left|\frac{f_1(x_0)}{f_2^2(x_0)}\right|\kappa_2^g(x_0)\Delta^{N_2^g}\\
			=& \left|\frac{1}{f_2(x_0)}\right|\kappa_1^g(x_0)\Delta^{N_1^g} + \left|\frac{f_1(x_0)}{f_2^2(x_0)}\right|\kappa_2^g(x_0)\Delta^{N_2^g}\\
			\le& K_F^g(x_0)\Delta^{N_F^g}.	
		\end{align*}
		That is, $\{\widetilde{F}_\Delta\}_{(0,\bar{\Delta}]}$ provides order-$N_F^g$ gradient accuracy of $F$ at $x_0$.
		
		$\hfill\qed$
	\end{proof}
	
    Theorem \ref{thm:qga} shows that, in terms of gradient accuracy at $x_0$, $\{\widetilde{F}_\Delta\}_{(0,\bar{\Delta}]}$ behaves at least as well as the worst behaved sub-model.
    
    Since the order-$\infty$ function accuracy assumption does not hold without loss of generality for the ‘near $x_0$’ case, the result is weaker.
	
	\begin{theorem}{(Gradient accuracy near $x_0$)}\label{thm:qgn}
		Let $f_1,f_2$ be in $\mathcal{C}^1$. Suppose for all $i\in\{1,2\}$, $\{\widetilde{f}_{i,\Delta}\}_{(0,\bar{\Delta}]}$ provides order-$N_i^f$ function accuracy of $f_i$ near $x_0$, and order-$N_i^g$ gradient accuracy of $f_i$ near $x_0$. If $N_2^f > 0$, then $\{\widetilde{F}_\Delta\}_{(0,\bar{\Delta}]}$ provides order-$N_F^{fg}$ gradient accuracy of $F$ near $x_0$, where $N_F^{fg} = \min\{N_1^f,N_2^f,N_1^g,N_2^g\}$.
	\end{theorem}
	
	\begin{proof}
	    Let $E_F^g = ||\nabla F(x)-\nabla\widetilde{F}_\Delta(x)||$. Changing $x_0$ to $x$ in the proof of Theorem \ref{thm:qga}, we obtain
	    \begin{align*}
	        E_F^g \le& \left\|\frac{\nabla\widetilde{f}_{2,\Delta}(x)}{f_2^2(x)\widetilde{f}_{2,\Delta}^2(x)}\right\|\left[\left|f_1(x)\left(f_2(x) + \widetilde{f}_{2,\Delta}(x)\right)\right|\kappa_2^f(x_0)\Delta^{N_2^f} + |f_2^2(x)|\kappa_1^f(x_0)\Delta^{N_1^f}\right]\\
			&+ \left\|\frac{\nabla \widetilde{f}_{1,\Delta}(x)}{f_2(x)\widetilde{f}_{2,\Delta}(x)}\right\|\kappa_2^f(x_0)\Delta^{N_2^f} + \left|\frac{1}{f_2(x)}\right|\kappa_1^g(x_0)\Delta^{N_1^g} + \left|\frac{f_1(x)}{f_2^2(x)}\right|\kappa_2^g(x_0)\Delta^{N_2^g}.
	    \end{align*}
	    Since all the coefficients above are bounded, the proof is complete.
	    
	    $\hfill\qed$
	\end{proof}

    Theorem \ref{thm:qgn} shows that, in terms of gradient accuracy near $x_0$, $\{\widetilde{F}_\Delta\}_{(0,\bar{\Delta}]}$ behaves at least as well as the worst behaved sub-model, where both function and gradient behaviors are considered.

	\begin{theorem}{(Hessian accuracy at $x_0$)}\label{thm:qha}
		Let $f_1,f_2$ be in $\mathcal{C}^2$. Suppose for all $i\in\{1,2\}$, $\{\widetilde{f}_{i,\Delta}\}_{(0,\bar{\Delta}]}$ provides order-$\infty$ function accuracy of $f_i$ at $x_0$ with constant $\kappa_i^f(x_0)=0$, order-$N_i^g$ gradient accuracy of $f_i$ at $x_0$ with constant $\kappa_i^g(x_0)$, and order-$N_i^h$ Hessian accuracy of $f_i$ at $x_0$ with constant $\kappa_i^h(x_0)$. Let $N_F^{gh}=\min\{N_1^g,N_2^g,N_1^h,N_2^h\}$. Then
		\begin{align*}
			||\nabla^2 F(x_0) - \nabla^2\widetilde{F}_\Delta(x_0)|| \le K_F^h(x_0)\Delta^{N_F^{gh}},
		\end{align*}
		where
		\begin{align*}
			K_F^h(x_0) =& \left| \frac{1}{f_2^2(x_0)}\right| (M_{\nabla f_2}(x_0) + M_{\nabla\widetilde{f}_2}(x_0))\kappa_1^g(x_0)\bar{\Delta}^{N_1^g - N_F^{gh}} + \left[ \left| \frac{1}{f_2^2(x_0)} \right| (M_{\nabla f_1}(x_0) + M_{\nabla\widetilde{f}_1}(x_0)) \right.\\
			&\left.+ \left| \frac{2f_1(x_0)}{f_2^3(x_0)} \right| (M_{\nabla f_2}(x_0) + M_{\nabla\widetilde{f}_2}(x_0))\right]\kappa_2^g(x_0)\bar{\Delta}^{N_2^g - N_F^{gh}} + \left| \frac{1}{f_2(x_0)} \right| \kappa_1^h(x_0)\bar{\Delta}^{N_1^h - N_F^{gh}}\\
			&+ \left| \frac{f_1(x_0)}{f_2^2(x_0)} \right| \kappa_2^h(x_0)\bar{\Delta}^{N_2^h - N_F^{gh}}.
		\end{align*}
		In particular, $\{\widetilde{F}_\Delta\}_{(0,\bar{\Delta}]}$ provides order-$N_F^{gh}$ Hessian accuracy for $F$ at $x_0$.
	\end{theorem}

	\begin{proof}
		Let $E_F^h = ||\nabla^2 F(x_0) - \nabla^2\widetilde{F}_\Delta(x_0)||$. Then
		\begin{align*}
			E_F^h =& \left\|\left(\frac{2f_1}{f_2^3}\nabla f_2\nabla f_2^\top - \frac{f_1}{f_2^2}\nabla^2 f_2 - \frac{1}{f_2^2}\nabla f_1\nabla f_2^\top - \frac{1}{f_2^2}\nabla f_2\nabla f_1^\top + \frac{1}{f_2}\nabla^2 f_1 - \frac{2\widetilde{f}_{1,\Delta}}{\widetilde{f}_{2,\Delta}^3}\nabla\widetilde{f}_{2,\Delta}\nabla\widetilde{f}_{2,\Delta}^\top\right.\right.\\
			&\left.\left. + \frac{\widetilde{f}_{1,\Delta}}{\widetilde{f}_{2,\Delta}^2}\nabla^2\widetilde{f}_{2,\Delta} + \frac{1}{\widetilde{f}_{2,\Delta}^2}\nabla \widetilde{f}_{1,\Delta}\nabla\widetilde{f}_{2,\Delta}^\top + \frac{1}{\widetilde{f}_{2,\Delta}^2}\nabla\widetilde{f}_{2,\Delta}\nabla\widetilde{f}_{1,\Delta}^\top - \frac{1}{\widetilde{f}_{2,\Delta}}\nabla^2\widetilde{f}_{1,\Delta}\right)(x_0)\right\|\\
			=& \left\|\left(\frac{2f_1}{f_2^3}\nabla f_2\nabla f_2^\top - \frac{2f_1}{f_2^3}\nabla\widetilde{f}_{2,\Delta}\nabla f_2^\top - \frac{f_1}{f_2^2}\nabla^2 f_2 + \frac{f_1}{f_2^2}\nabla^2\widetilde{f}_{2,\Delta} - \frac{1}{f_2^2}\nabla f_1\nabla f_2^\top + \frac{1}{f_2^2}\nabla\widetilde{f}_{1,\Delta}\nabla f_2^\top\right.\right.\\ 
			&\left.\left. - \frac{1}{f_2^2}\nabla f_2\nabla f_1^\top + \frac{1}{f_2^2}\nabla\widetilde{f}_{2,\Delta}\nabla f_1^\top + \frac{1}{f_2}\nabla^2 f_1 - \frac{1}{f_2}\nabla^2\widetilde{f}_{1,\Delta} -\frac{2\widetilde{f}_{1,\Delta}}{\widetilde{f}_{2,\Delta}^3}\nabla\widetilde{f}_{2,\Delta}\nabla\widetilde{f}_{2,\Delta}^\top + \frac{2f_1}{f_2^3}\nabla\widetilde{f}_{2,\Delta}\nabla f_2^\top\right.\right.\\
			&\left.\left. + \frac{\widetilde{f}_{1,\Delta}}{\widetilde{f}_{2,\Delta}^2}\nabla^2\widetilde{f}_{2,\Delta} - \frac{f_1}{f_2^2}\nabla^2\widetilde{f}_{2,\Delta} + \frac{1}{\widetilde{f}_{2,\Delta}^2}\nabla \widetilde{f}_{1,\Delta}\nabla\widetilde{f}_{2,\Delta}^\top - \frac{1}{f_2^2}\nabla\widetilde{f}_{1,\Delta}\nabla f_2^\top + \frac{1}{\widetilde{f}_{2,\Delta}^2}\nabla\widetilde{f}_{2,\Delta}\nabla\widetilde{f}_{1,\Delta}^\top\right.\right.\\
			&\left.\left. - \frac{1}{f_2^2}\nabla\widetilde{f}_{2,\Delta}\nabla f_1^\top - \frac{1}{\widetilde{f}_{2,\Delta}}\nabla^2\widetilde{f}_{1,\Delta} + \frac{1}{f_2}\nabla^2\widetilde{f}_{1,\Delta}\right)(x_0)\right\|\\
			=& \left\|\left(\frac{\widetilde{f}_{1,\Delta}}{\widetilde{f}_{2,\Delta}^2}\nabla^2\widetilde{f}_{2,\Delta} - \frac{f_1}{f_2^2}\nabla^2\widetilde{f}_{2,\Delta} + \frac{1}{f_2}\nabla^2\widetilde{f}_{1,\Delta} - \frac{1}{\widetilde{f}_{2,\Delta}}\nabla^2\widetilde{f}_{1,\Delta} + \frac{1}{\widetilde{f}_{2,\Delta}^2}\nabla\widetilde{f}_{2,\Delta}\nabla\widetilde{f}_{1,\Delta}^\top - \frac{1}{f_2^2}\nabla\widetilde{f}_{2,\Delta}\nabla f_1^\top\right.\right.\\
			&\left.\left. + \frac{2f_1}{f_2^3}\nabla\widetilde{f}_{2,\Delta}\nabla f_2^\top - \frac{2\widetilde{f}_{1,\Delta}}{\widetilde{f}_{2,\Delta}^3}\nabla\widetilde{f}_{2,\Delta}\nabla\widetilde{f}_{2,\Delta}^\top + \frac{1}{\widetilde{f}_{2,\Delta}^2}\nabla \widetilde{f}_{1,\Delta}\nabla\widetilde{f}_{2,\Delta}^\top - \frac{1}{f_2^2}\nabla\widetilde{f}_{1,\Delta}\nabla f_2^\top + \frac{1}{f_2^2}\nabla\widetilde{f}_{1,\Delta}\nabla f_2^\top\right.\right.\\
			&\left.\left. - \frac{1}{f_2^2}\nabla f_1\nabla f_2^\top + \frac{1}{f_2^2}\nabla\widetilde{f}_{2,\Delta}\nabla f_1^\top - \frac{1}{f_2^2}\nabla f_2\nabla f_1^\top + \frac{2f_1}{f_2^3}\nabla f_2\nabla f_2^\top - \frac{2f_1}{f_2^3}\nabla\widetilde{f}_{2,\Delta}\nabla f_2^\top + \frac{1}{f_2}\nabla^2 f_1\right.\right.\\
			&\left.\left. - \frac{1}{f_2}\nabla^2\widetilde{f}_{1,\Delta} + \frac{f_1}{f_2^2}\nabla^2\widetilde{f}_{2,\Delta} - \frac{f_1}{f_2^2}\nabla^2 f_2\right)(x_0)\right\|\\
			\le& \left(\left\| \frac{\nabla^2 \widetilde{f}_{2,\Delta}(x_0)}{f_2^2(x_0)} \right\|+ \left\| \frac{2\nabla\widetilde{f}_{2,\Delta}(x_0)\nabla\widetilde{f}_{2,\Delta}(x_0)^\top}{f_2^3(x_0)} \right\|\right) \kappa_1^f(x_0)\Delta^{N_1^f} + \left[ \left\| \frac{\nabla^2 \widetilde{f}_{1,\Delta}(x_0)}{f_2(x_0)\widetilde{f}_{2,\Delta}(x_0)} \right\|\right.\\
			&\left.+ \left\| \frac{\widetilde{f}_{1,\Delta}(f_2 + \widetilde{f}_{2,\Delta})\nabla^2 \widetilde{f}_{2,\Delta}}{f_2^2\widetilde{f}_{2,\Delta}^2}(x_0) \right\| + \left\| \frac{(f_2 + \widetilde{f}_{2,\Delta})(\nabla\widetilde{f}_{1,\Delta}\nabla\widetilde{f}_{2,\Delta}^\top + \nabla\widetilde{f}_{2,\Delta}\nabla\widetilde{f}_{1,\Delta}^\top)}{f_2^2\widetilde{f}_{2,\Delta}^2}(x_0) \right\|\right.\\
			&\left.+ \left\| \frac{2\widetilde{f}_{1,\Delta}(f_2^2 + f_2\widetilde{f}_{2,\Delta} + \widetilde{f}_{2,\Delta}^2)\nabla\widetilde{f}_{2,\Delta}\nabla\widetilde{f}_{2,\Delta}^\top}{f_2^3\widetilde{f}_{2,\Delta}^3}(x_0) \right\|\right] \kappa_2^f(x_0)\Delta^{N_2^f} + \left| \frac{1}{f_2^2(x_0)}\right| (||\nabla f_2(x_0)||\\
			& + ||\nabla\widetilde{f}_{2,\Delta}(x_0)||)\kappa_1^g(x_0)\Delta^{N_1^g} + \left[ \left| \frac{1}{f_2^2(x_0)} \right| (||\nabla f_1(x_0)|| + ||\nabla\widetilde{f}_{1,\Delta}(x_0)||) + \left| \frac{2f_1(x_0)}{f_2^3(x_0)} \right| (||\nabla f_2(x_0)||\right.\\
			&\left. + ||\nabla\widetilde{f}_{2,\Delta}(x_0)||)\right]\kappa_2^g(x_0)\Delta^{N_2^g} + \left| \frac{1}{f_2(x_0)} \right| \kappa_1^h(x_0)\Delta^{N_1^h} + \left| \frac{f_1(x_0)}{f_2^2(x_0)} \right| \kappa_2^h(x_0)\Delta^{N_2^h}\\
			=& \left| \frac{1}{f_2^2(x_0)}\right| (||\nabla f_2(x_0)|| + ||\nabla\widetilde{f}_{2,\Delta}(x_0)||)\kappa_1^g(x_0)\Delta^{N_1^g} + \left[ \left| \frac{1}{f_2^2(x_0)} \right| (||\nabla f_1(x_0)|| + ||\nabla\widetilde{f}_{1,\Delta}(x_0)||)\right.\\
			&\left.+ \left| \frac{2f_1(x_0)}{f_2^3(x_0)} \right| (||\nabla f_2(x_0)|| + ||\nabla\widetilde{f}_{2,\Delta}(x_0)||)\right]\kappa_2^g(x_0)\Delta^{N_2^g} + \left| \frac{1}{f_2(x_0)} \right| \kappa_1^h(x_0)\Delta^{N_1^h}\\
			&+ \left| \frac{f_1(x_0)}{f_2^2(x_0)} \right| \kappa_2^h(x_0)\Delta^{N_2^h}\\
			\le& K_F^h(x_0)\Delta^{N_F^{gh}}.
	    \end{align*}
		
	    $\hfill\qed$
	\end{proof}

    Theorem \ref{thm:qha} shows that, in terms of Hessian accuracy at $x_0$, $\{\widetilde{F}_\Delta\}_{(0,\bar{\Delta}]}$ behaves at least as well as the worst behaved sub-model, where both gradient and Hessian behaviors are considered.
	
	\begin{theorem}{(Hessian accuracy near $x_0$)}\label{thm:qhn}
		Let $f_1,f_2$ be in $\mathcal{C}^2$. Suppose for all $i\in\{1,2\}$, $\{\widetilde{f}_{i,\Delta}\}_{(0,\bar{\Delta}]}$ provides order-$N_i^f$ function accuracy of $f_i$ near $x_0$, order-$N_i^g$ gradient accuracy of $f_i$ near $x_0$, and order-$N_i^h$ Hessian accuracy of $f_i$ near $x_0$. If $N_2^f > 0$, then $\{\widetilde{F}_\Delta\}_{(0,\bar{\Delta}]}$ provides order-$N_F^{fgh}$ Hessian accuracy for $F$ near $x_0$, where $N_F^{fgh}=\min\{N_1^f,N_2^f,N_1^g,N_2^g,N_1^h,N_2^h\}$.
	\end{theorem}
	
	\begin{proof}
	    Let $E_F^h = ||\nabla^2 F(x)-\nabla^2\widetilde{F}_\Delta(x)||$. Changing $x_0$ to $x$ in the proof of Theorem \ref{thm:qha}, we obtain
	    \begin{align*}
	        E_F^h \le& \left(\left\| \frac{\nabla^2 \widetilde{f}_{2,\Delta}(x)}{f_2^2(x)} \right\|+ \left\| \frac{2\nabla\widetilde{f}_{2,\Delta}(x)\nabla\widetilde{f}_{2,\Delta}(x)^\top}{f_2^3(x)} \right\|\right) \kappa_1^f(x_0)\Delta^{N_1^f} + \left[ \left\| \frac{\nabla^2 \widetilde{f}_{1,\Delta}(x)}{f_2(x)\widetilde{f}_{2,\Delta}(x)} \right\|\right.\\
			&\left.+ \left\| \frac{\widetilde{f}_{1,\Delta}(f_2 + \widetilde{f}_{2,\Delta})\nabla^2 \widetilde{f}_{2,\Delta}}{f_2^2\widetilde{f}_{2,\Delta}^2}(x) \right\| + \left\| \frac{(f_2 + \widetilde{f}_{2,\Delta})(\nabla\widetilde{f}_{1,\Delta}\nabla\widetilde{f}_{2,\Delta}^\top + \nabla\widetilde{f}_{2,\Delta}\nabla\widetilde{f}_{1,\Delta}^\top)}{f_2^2\widetilde{f}_{2,\Delta}^2}(x) \right\|\right.\\
			&\left.+ \left\| \frac{2\widetilde{f}_{1,\Delta}(f_2^2 + f_2\widetilde{f}_{2,\Delta} + \widetilde{f}_{2,\Delta}^2)\nabla\widetilde{f}_{2,\Delta}\nabla\widetilde{f}_{2,\Delta}^\top}{f_2^3\widetilde{f}_{2,\Delta}^3}(x) \right\|\right] \kappa_2^f(x_0)\Delta^{N_2^f} + \left| \frac{1}{f_2^2(x)}\right| (||\nabla f_2(x)||\\
			& + ||\nabla\widetilde{f}_{2,\Delta}(x)||)\kappa_1^g(x_0)\Delta^{N_1^g} + \left[ \left| \frac{1}{f_2^2(x)} \right| (||\nabla f_1(x)|| + ||\nabla\widetilde{f}_{1,\Delta}(x)||) + \left| \frac{2f_1(x)}{f_2^3(x)} \right| (||\nabla f_2(x)||\right.\\
			&\left. + ||\nabla\widetilde{f}_{2,\Delta}(x)||)\right]\kappa_2^g(x_0)\Delta^{N_2^g} + \left| \frac{1}{f_2(x)} \right| \kappa_1^h(x_0)\Delta^{N_1^h} + \left| \frac{f_1(x)}{f_2^2(x)} \right| \kappa_2^h(x_0)\Delta^{N_2^h}.
	    \end{align*}
	    Since all the coefficients above are bounded, the proof is complete.
	    
	    $\hfill\qed$
	\end{proof}

    Theorem \ref{thm:qhn} shows that, in terms of Hessian accuracy near $x_0$, $\{\widetilde{F}_\Delta\}_{(0,\bar{\Delta}]}$ behaves at least as well as the worst behaved sub-model, where all function, gradient, and Hessian behaviors are considered.

	The next example shows the necessity of $N_2^f > 0$ in this section. In particular, if $N_2^f = 0$, then the function error of $\{\widetilde{F}_\Delta\}_{(0,\bar{\Delta}]}$ can be infinitely large.

	\begin{example}
		Suppose $\{\widetilde{f}_{1,\Delta}\}_{(0,\bar{\Delta}]}$ provides order-$N_1^f$ function accuracy of $f$ at $x_0$, and $\{\widetilde{f}_{2,\Delta}\}_{(0,\bar{\Delta}]}$ provide order-$N_2^f$ function accuracy of $f_2$ at $x_0$. If $N_2^f = 0$, then $|F(x_0) - \widetilde{F}_\Delta(x_0)|$ can be infinitely large.
	\end{example}

	\begin{proof}
		Let $f_1(x) = 1, \widetilde{f}_{1,\Delta}(x) = 1, f_2(x)=1, \widetilde{f}_{2,\Delta}(x) = \Delta$. Then, regardless of $N_1^f$, we have
		\begin{align*}
			|f_1(x_0) - \widetilde{f}_{1,\Delta}(x_0)| =& 0 \le \Delta^{N_1^f},\\
			|F(x_0) - \widetilde{F}_\Delta(x_0)| =& |1 - \frac{1}{\Delta}|.
		\end{align*}
		Since $\lim\limits_{\Delta\to0} |1 - \frac{1}{\Delta}| =\infty$, the proof is complete.
		
		$\hfill\qed$
	\end{proof}

	Notice that $\{\widetilde{f}_{1,\Delta}\}_{(0,\bar{\Delta}]}$ and $\{\widetilde{f}_{2,\Delta}\}_{(0,\bar{\Delta}]}$ in the example above also obtain the true gradient and Hessian of $f_1$ and $f_2$. Therefore, even if we have more information in some theorems in this section, we still need $N_2^f > 0$.

	\section{Error Analysis of the Composition of Model Functions}
	
	In this section, we consider the case where the objective function $F$ can be expressed as the composition of a component mapping and a component function. Throughout this section, we assume that $f_1$ is a mapping from $\mathbb{R}^d$ to $\mathbb{R}^m$, and $f_2$ is a function from $\mathbb{R}^m$ to $\mathbb{R}$. Let $F = f_2\circ f_1,\widetilde{F}_\Delta = \widetilde{f}_{2,\Delta}\circ\widetilde{f}_{1,\Delta}$. For the sake of clarity, we denote $\nabla g(f(\cdot))$ by $\nabla g\big|_f(\cdot)$ in this section. The next two lemmas are crucial in the proofs of this section.
	
	\begin{lemma}\label{lem:lip}
		If a function $f\in\mathcal{C}^1$, then it is Lipschitz continuous on any convex compact set $K$ with constant $\max\limits_{x\in K}||\nabla f(x)||$.
	\end{lemma}
	
	\begin{proof}
		See \cite[Thm 9.7]{rockafellar2009variational}.
		
		$\hfill\qed$
	\end{proof}

	\begin{lemma}\label{lem:img}
		Suppose the mapping $g = (g_1,...,g_m)^\top\in\mathcal{C}^1$. Let $g(B_\Delta(x_0))$ be the image of $B_\Delta(x_0)$ under $g$. Let $M_{\nabla g_i}(x_0) = \max\limits_{x\in\overline{B_{\bar{\Delta}}(x_0)}} ||\nabla g_i(x)||,\overline{M}_{\nabla g}(x_0) = ||(M_{\nabla g_1}(x_0),...,M_{\nabla g_n}(x_0))^\top||$. Then for all $\Delta\in(0,\bar{\Delta}]$, $g(B_\Delta(x_0))\subseteq B_{\overline{M}_{\nabla g}(x_0)\Delta}(g(x_0))$.
	\end{lemma}
	
	\begin{proof}
		For all $i\in\{1,...,m\}$ and $x\in B_\Delta(x_0)$ we have $$|g_i(x) - g_i(x_0)| \le M_{\nabla g_i}(x_0) ||x - x_0||.$$ Therefore $$||g(x) - g(x_0)|| \le \overline{M}_{\nabla g}(x_0) ||x - x_0||.$$
		That is, $g(B_\Delta(x_0))\subseteq B_{\overline{M}_{\nabla g}(x_0)\Delta}(g(x_0))$.
		
		$\hfill\qed$
	\end{proof}

	Using the two lemmas above, we develop the overall function error bound.
	
	\begin{theorem}{(Function accuracy at/near $x_0$)}\label{thm:cf}
		Let $f_1$ be in $\mathcal{C}^1$. Suppose $\{\widetilde{f}_{1,\Delta}\}_{(0,\bar{\Delta}]}$ provides order-$N_1^f$ function accuracy of $f_1$ at/near $x_0$ with constant $\kappa_1^f(x_0)$, and $\{\widetilde{f}_{2,\Delta}\}_{(0,\bar{\Delta}]}$ provides order-$N_2^f$ function accuracy of $f_2$ at/near $f_1(x_0)$ with constant $\kappa_2^f(f_1(x_0))$. Let $N_F^f = \min\{N_1^f,N_2^f\}$. If $N_1^f > 0$, $\{\widetilde{f}_{2,\Delta}\}_{(0,\bar{\Delta}]}\subseteq\mathcal{C}^1$, and $\{||\nabla\widetilde{f}_{2,\Delta}||\}_{(0,\bar{\Delta}]}$ is uniformly bounded at/near $f_1(x_0)$ by $M_{\nabla\widetilde{f}_2}(f_1(x_0))$. Then
		
		$\ $(i) if the accuracy provided is at $f_1(x_0)$ and $x_0$, then
		\begin{align*}
			|F(x_0) - \widetilde{F}_\Delta(x_0)| \le K_F^f(x_0)\Delta^{N_F^f};
		\end{align*}
		
		(ii) if the accuracy provided is near $f_1(x_0)$ and $x_0$, then for all $x\in B_\Delta(x_0)$
		\begin{align*}
			|F(x) - \widetilde{F}_\Delta(x)| \le K_F^f(x_0)\Delta^{N_F^f};
		\end{align*}
		where
		\begin{align*}
			K_F^f(x_0) = M_{\nabla\widetilde{f}_2}(f_1(x_0))\kappa_1^f(x_0)\bar{\Delta}^{N_1^f-N_F^f} + \overline{M}_{\nabla f_1}(x_0)^{N_2^f}\kappa_2^f(f_1(x_0))\bar{\Delta}^{N_2^f-N_F^f}.
		\end{align*}
		In particular, $\{\widetilde{F}_\Delta\}_{(0,\bar{\Delta}]}$ provides order-$N_F^f$ function accuracy of $F$ at/near $x_0$.
	\end{theorem}
	
	\begin{proof}
		We only prove case (ii), since the proof of case (i) can easily be obtained by changing all $x$ to $x_0$. Let $$S = f_1(B_{\bar{\Delta}}(x_0))\cup\bigcup\limits_{\Delta\in(0,\bar{\Delta}]}\widetilde{f}_{1,\Delta}(B_{\bar{\Delta}}(x_0)),$$ where $f_1(B_{\bar{\Delta}}(x_0))$ and $\widetilde{f}_{1,\Delta}(B_{\bar{\Delta}}(x_0))$ are the image of $B_{\bar{\Delta}}(x_0)$ under $f_1$ and $\widetilde{f}_{1,\Delta}$. Let $K$ be the closed convex hull of $S$. Since $\{\widetilde{f}_{1,\Delta}\}_{(0,\bar{\Delta}]}$ is uniformly bounded, it is easy to see that $K$ is a convex compact set. Since $N_1^f > 0$, by Lemma \ref{lem:img}, we can choose a sufficiently small $\bar{\Delta}$ such that $\{\widetilde{f}_{2,\Delta}\}_{(0,\bar{\Delta}]}$ is uniformly bounded in $K$. By Lemma \ref{lem:lip}, we get for all $\Delta\in(0,\bar{\Delta}],\widetilde{f}_{2,\Delta}$ is Lipschitz continuous on $K$ with constant $M_{\nabla\widetilde{f}_2}(f_1(x_0))$. Therefore
		\begin{align*}
			|F(x) - \widetilde{F}_\Delta(x)| =& |(f_2\circ f_1 - \widetilde{f}_{2,\Delta}\circ\widetilde{f}_{1,\Delta})(x)|\\
			=& |(f_2\circ f_1 - \widetilde{f}_{2,\Delta}\circ f_1 + \widetilde{f}_{2,\Delta}\circ f_1 - \widetilde{f}_{2,\Delta}\circ\widetilde{f}_{1,\Delta})(x)|\\
			=& |(\widetilde{f}_{2,\Delta}\circ f_1 - \widetilde{f}_{2,\Delta}\circ\widetilde{f}_{1,\Delta} + f_2\circ f_1 - \widetilde{f}_{2,\Delta}\circ f_1)(x)|\\
			\le& M_{\nabla\widetilde{f}_2}(f_1(x_0))\kappa_1^f(x_0)\Delta^{N_1^f} + \overline{M}_{\nabla f_1}(x_0)^{N_2^f}\kappa_2^f(f_1(x_0))\Delta^{N_2^f}\\
			\le& K_F^f(x_0)\Delta^{N_F^f}.
		\end{align*}
		That is, $\{\widetilde{F}_\Delta\}_{(0,\bar{\Delta}]}$ provides order-$N_F^f$ function accuracy of $F$ at/near $x_0$.
		
		$\hfill\qed$
	\end{proof}
	
	Theorem \ref{thm:cf} shows that, in terms of function accuracy, $\{\widetilde{F}_\Delta\}_{(0,\bar{\Delta}]}$ behaves at least as well as the worst behaved sub-model.
	
	Also, we have error bounds for gradient and Hessian.
	
	\begin{theorem}{(Gradient accuracy at $x_0$)}\label{thm:cga}
		Let $f_1,f_2$ be in $\mathcal{C}^1$. Suppose $\{\widetilde{f}_{1,\Delta}\}_{(0,\bar{\Delta}]}$ provides order-$\infty$ function accuracy of $f_1$ at $x_0$ with constant $\kappa_1^f(x_0)=0$, and order-$N_1^g$ gradient accuracy of $f_1$ at $x_0$ with constant $\kappa_1^g(x_0)$. Suppose $\{\widetilde{f}_{2,\Delta}\}_{(0,\bar{\Delta}]}$ provides order-$\infty$ function accuracy of $f_2$ at $f_1(x_0)$ with constant $\kappa_2^f(f_1(x_0))=0$, and order-$N_2^g$ gradient accuracy of $f_2$ at $f_1(x_0)$ with constant $\kappa_2^g(f_1(x_0))$. Let $N_F^g=\min\{N_1^g,N_2^g\}$. Then
		\begin{align*}
			||\nabla F(x_0) - \nabla\widetilde{F}_\Delta(x_0)|| \le K_F^g(x_0)\Delta^{N_F^g},
		\end{align*}
		where $$K_F^g(x_0) = ||\nabla f_2(f_1(x_0))||\kappa_1^g(x_0)\bar{\Delta}^{N_1^g - N_F^g} + M_{\nabla\widetilde{f}_1}(x_0)\kappa_2^g(f_1(x_0))\bar{\Delta}^{N_2^g - N_F^g}.$$
		In particular, $\{\widetilde{F}_\Delta\}_{(0,\bar{\Delta}]}$ provides order-$N_F^g$ gradient accuracy of $F$ at $x_0$.
	\end{theorem}
	
	\begin{proof}
		Let $E_F^g = ||\nabla F(x_0) - \nabla\widetilde{F}(x_0)||$. Then
		\begin{align*}
			E_F^g =& ||(\nabla f_1^\top\nabla f_2\big|_{f_1} - \nabla {\widetilde{f}_{1,\Delta}}^\top\nabla \widetilde{f}_{2,\Delta}\big|_{\widetilde{f}_{1,\Delta}})(x_0)||\\
			=& ||(\nabla f_1^\top\nabla f_2\big|_{f_1} - \nabla {\widetilde{f}_{1,\Delta}}^\top\nabla f_2\big|_{f_1} + \nabla {\widetilde{f}_{1,\Delta}}^\top\nabla f_2\big|_{f_1} - \nabla {\widetilde{f}_{1,\Delta}}^\top\nabla \widetilde{f}_{2,\Delta}\big|_{f_1} + \nabla {\widetilde{f}_{1,\Delta}}^\top\nabla \widetilde{f}_{2,\Delta}\big|_{f_1}\\
			&- \nabla {\widetilde{f}_{1,\Delta}}^\top\nabla \widetilde{f}_{2,\Delta}\big|_{\widetilde{f}_{1,\Delta}})(x_0)||\\
			=& ||(\nabla f_1^\top\nabla f_2\big|_{f_1} - \nabla {\widetilde{f}_{1,\Delta}}^\top\nabla f_2\big|_{f_1} + \nabla {\widetilde{f}_{1,\Delta}}^\top\nabla f_2\big|_{f_1} -  \nabla {\widetilde{f}_{1,\Delta}}^\top\nabla \widetilde{f}_{2,\Delta}\big|_{f_1})(x_0)||\\
			\le& ||\nabla f_2(f_1(x_0))||\kappa_1^g(x_0)\Delta^{N_1^g} + ||\nabla {\widetilde{f}_{1,\Delta}}(x_0)||\kappa_2^g(f_1(x_0))\Delta^{N_2^g}\\
			\le& K_F^g(x_0)\Delta^{N_F^g}.
		\end{align*}
		That is, $\{\widetilde{F}_\Delta\}_{(0,\bar{\Delta}]}$ provides order-$N_F^g$ gradient accuracy of $F$ at $x_0$.
		
		$\hfill\qed$
	\end{proof}
	
	Theorem \ref{thm:cga} shows that, in terms of gradient accuracy at $x_0$, $\{\widetilde{F}_\Delta\}_{(0,\bar{\Delta}]}$ behaves at least as well as the worst behaved sub-model.
	
	Since the order-$\infty$ function accuracy assumption does not hold without loss of generality for the ‘near $x_0$’ case, the result is weaker.
	
	\begin{theorem}{(Gradient accuracy near $x_0$)}\label{thm:cgn}
		Let $f_1,f_2$ be in $\mathcal{C}^1$. Suppose $\{\widetilde{f}_{1,\Delta}\}_{(0,\bar{\Delta}]}$ provides order-$N_1^f$ function accuracy of $f_1$ near $x_0$, and order-$N_1^g$ gradient accuracy of $f_1$ near $x_0$. Suppose $\{\widetilde{f}_{2,\Delta}\}_{(0,\bar{\Delta}]}$ provides order-$N_2^f$ function accuracy of $f_2$ near $f_1(x_0)$, and order-$N_2^g$ gradient accuracy of $f_2$ near $f_1(x_0)$. If $N_1^f > 0$, $\{\widetilde{f}_{2,\Delta}\}_{(0,\bar{\Delta}]}\subseteq\mathcal{C}^2$, and $\{||\nabla^2\widetilde{f}_{2,\Delta}||\}_{(0,\bar{\Delta}]}$ is uniformly bounded near $f_1(x_0)$ by $M_{\nabla^2\widetilde{f}_2}(f_1(x_0))$, then $\{\widetilde{F}_\Delta\}_{(0,\bar{\Delta}]}$ provides order-$N_F^{fg}$ gradient accuracy of $F$ near $x_0$, where $N_F^{fg}=\min\{N_1^f,N_2^f,N_1^g,N_2^g\}$.
	\end{theorem}
	
	\begin{proof}
		Let $E_F^g = ||\nabla F(x) - \nabla\widetilde{F}(x)||$. Changing $x_0$ to $x$ in the proof of Theorem \ref{thm:cga}, we obtain
		\begin{align*}
		    E_F^g =& ||(\nabla f_1^\top\nabla f_2\big|_{f_1} - \nabla {\widetilde{f}_{1,\Delta}}^\top\nabla f_2\big|_{f_1} + \nabla {\widetilde{f}_{1,\Delta}}^\top\nabla f_2\big|_{f_1} - \nabla {\widetilde{f}_{1,\Delta}}^\top\nabla \widetilde{f}_{2,\Delta}\big|_{f_1} + \nabla {\widetilde{f}_{1,\Delta}}^\top\nabla \widetilde{f}_{2,\Delta}\big|_{f_1}\\
			&- \nabla {\widetilde{f}_{1,\Delta}}^\top\nabla \widetilde{f}_{2,\Delta}\big|_{\widetilde{f}_{1,\Delta}})(x)||.
		\end{align*}
		Noticing that $$|| \nabla f_2(f_1(x)) -  \nabla \widetilde{f}_{2,\Delta}(f_1(x))|| \le \overline{M}_{\nabla f_1}(x_0)^{N_2^g}\kappa_2^g(f_1(x_0))\Delta^{N_2^g},$$ and $$||\nabla\widetilde{f}_{2,\Delta}(f_1(x)) - \nabla \widetilde{f}_{2,\Delta}(\widetilde{f}_{1,\Delta}(x))|| \le M_{\nabla^2\widetilde{f}_2}(f_1(x_0))\kappa_1^f(x_0)\Delta^{N_1^f},$$ we can get the conclusion.
		
		$\hfill\qed$
	\end{proof}
	
	Theorem \ref{thm:cgn} shows that, in terms of gradient accuracy near $x_0$, $\{\widetilde{F}_\Delta\}_{(0,\bar{\Delta}]}$ behaves at least as well as the worst behaved sub-model, where both function and gradient behaviors are considered.
	
	\begin{theorem}{(Hessian accuracy at $x_0$)}\label{thm:cha}
		Let $f_1,f_2$ be in $\mathcal{C}^2$. Suppose $\{\widetilde{f}_{1,\Delta}\}_{(0,\bar{\Delta}]}$ provides order-$\infty$ function accuracy of $f_1$ at $x_0$ with constant $\kappa_1^f(x_0)=0$, order-$N_1^g$ gradient accuracy of $f_1$ at $x_0$ with constant $\kappa_1^g(x_0)$, and order-$N_1^h$ Hessian accuracy of $f_1$ at $x_0$ with constant $\kappa_1^h(x_0)$. Suppose $\{\widetilde{f}_{2,\Delta}\}_{(0,\bar{\Delta}]}$ provides order-$\infty$ function accuracy of $f_2$ at $f_1(x_0)$ with constant $\kappa_2^f(f_1(x_0))=0$, order-$N_2^g$ gradient accuracy of $f_2$ at $f_1(x_0)$ with constant $\kappa_2^g(f_1(x_0))$, and order-$N_2^h$ Hessian accuracy of $f_2$ at $f_1(x_0)$ with constant $\kappa_2^h(f_1(x_0))$. Let $N_F^{gh}=\min\{N_1^g,N_2^g,N_1^h,N_2^h\}$. Then
		\begin{align*}
			||\nabla^2 F(x_0) - \nabla^2\widetilde{F}_\Delta(x_0)|| \le K_F^h(x_0)\Delta^{N_F^{gh}},
		\end{align*}
		where
		\begin{align*}
			K_F^h(x_0) =& M_{\nabla^2\widetilde{f}_2}(f_1(x_0))(||\nabla f_1(x_0)|| + M_{\nabla\widetilde{f}_1}(x_0))\kappa_1^g(x_0)\bar{\Delta}^{N_1^g - N_F^{gh}} + M_{\nabla^2\widetilde{f}_1}(x_0)\kappa_2^g(f_1(x_0))\bar{\Delta}^{N_2^g - N_F^{gh}}\\
			&+ (||\nabla f_2(f_1(x_0))|| + M_{\nabla\widetilde{f}_2}(f_1(x_0)))\kappa_1^h(x_0)\bar{\Delta}^{N_1^h - N_F^{gh}} + ||\nabla f_1(x_0)||^2\kappa_2^h(f_1(x_0))\bar{\Delta}^{N_2^h - N_F^{gh}}.
		\end{align*}
		In particular, $\{\widetilde{F}_\Delta\}_{(0,\bar{\Delta}]}$ provides order-$N_F^{gh}$ Hessian accuracy of $F$ at $x_0$.
	\end{theorem}
	
	\begin{proof}
		Let $E_F^h = ||\nabla^2 F(x_0) - \nabla^2\widetilde{F}_\Delta(x_0)||$. Then
		\begin{align*}
			E_F^h =& ||(\nabla f_1^\top\nabla^2f_2\big|_{f_1}\nabla f_1 + \nabla f_2\big|_{f_1}^\top\nabla^2 f_1 - \nabla\widetilde{f}_{1,\Delta}^\top\nabla^2\widetilde{f}_{2,\Delta}\big|_{\widetilde{f}_{1,\Delta}}\nabla\widetilde{f}_{1,\Delta} - \nabla\widetilde{f}_{2,\Delta}\big|_{\widetilde{f}_{1,\Delta}}^\top\nabla^2\widetilde{f}_{1,\Delta})(x_0)||\\
			=& ||(\nabla f_1^\top\nabla^2f_2\big|_{f_1}\nabla f_1 - \nabla f_1^\top\nabla^2\widetilde{f}_{2,\Delta}\big|_{f_1}\nabla f_1 + \nabla f_2\big|_{f_1}^\top\nabla^2 f_1 - \nabla f_2\big|_{f_1}^\top\nabla^2\widetilde{f}_{1,\Delta}\\
			& - \nabla\widetilde{f}_{1,\Delta}^\top\nabla^2\widetilde{f}_{2,\Delta}\big|_{\widetilde{f}_{1,\Delta}}\nabla\widetilde{f}_{1,\Delta} + \nabla f_1^\top\nabla^2\widetilde{f}_{2,\Delta}\big|_{\widetilde{f}_{1,\Delta}}\nabla\widetilde{f}_{1,\Delta} - \nabla\widetilde{f}_{2,\Delta}\big|_{\widetilde{f}_{1,\Delta}}^\top\nabla^2\widetilde{f}_{1,\Delta} + \nabla\widetilde{f}_{2,\Delta}\big|_{f_1}^\top\nabla^2\widetilde{f}_{1,\Delta}\\
			&+ \nabla f_1^\top\nabla^2\widetilde{f}_{2,\Delta}\big|_{f_1}\nabla f_1 - \nabla f_1^\top\nabla^2\widetilde{f}_{2,\Delta}\big|_{\widetilde{f}_{1,\Delta}}\nabla\widetilde{f}_{1,\Delta} + \nabla f_2\big|_{f_1}^\top\nabla^2\widetilde{f}_{1,\Delta} - \nabla\widetilde{f}_{2,\Delta}\big|_{f_1}^\top\nabla^2\widetilde{f}_{1,\Delta})(x_0)||\\
			=& ||(\nabla\widetilde{f}_{2,\Delta}\big|_{f_1}^\top\nabla^2\widetilde{f}_{1,\Delta} - \nabla\widetilde{f}_{2,\Delta}\big|_{\widetilde{f}_{1,\Delta}}^\top\nabla^2\widetilde{f}_{1,\Delta} + \nabla f_1^\top\nabla^2\widetilde{f}_{2,\Delta}\big|_{f_1}\nabla f_1 - \nabla f_1^\top\nabla^2\widetilde{f}_{2,\Delta}\big|_{\widetilde{f}_{1,\Delta}}\nabla\widetilde{f}_{1,\Delta}\\
			&+ \nabla f_1^\top\nabla^2\widetilde{f}_{2,\Delta}\big|_{\widetilde{f}_{1,\Delta}}\nabla\widetilde{f}_{1,\Delta} - \nabla\widetilde{f}_{1,\Delta}^\top\nabla^2\widetilde{f}_{2,\Delta}\big|_{\widetilde{f}_{1,\Delta}}\nabla\widetilde{f}_{1,\Delta} + \nabla f_2\big|_{f_1}^\top\nabla^2\widetilde{f}_{1,\Delta} - \nabla\widetilde{f}_{2,\Delta}\big|_{f_1}^\top\nabla^2\widetilde{f}_{1,\Delta}\\
			&+ \nabla f_2\big|_{f_1}^\top\nabla^2 f_1 - \nabla f_2\big|_{f_1}^\top\nabla^2\widetilde{f}_{1,\Delta} + \nabla f_1^\top\nabla^2f_2\big|_{f_1}\nabla f_1 - \nabla f_1^\top\nabla^2\widetilde{f}_{2,\Delta}\big|_{f_1}\nabla f_1)(x_0)||\\
			\le& ||\nabla^2\widetilde{f}_{2,\Delta}(f_1(x_0))||(||\nabla f_1(x_0)|| + ||\nabla\widetilde{f}_{1,\Delta}(x_0)||)\kappa_1^g(x_0)\Delta^{N_1^g} + ||\nabla^2\widetilde{f}_{1,\Delta}(x_0)||\kappa_2^g(f_1(x_0))\Delta^{N_2^g}\\
			&+ ||\nabla f_2(f_1(x_0))||\kappa_1^h(x_0)\Delta^{N_1^h} + ||\nabla f_1(x_0)||^2\kappa_2^h(f_1(x_0))\Delta^{N_2^h}\\
			\le& K_F^h(x_0)\Delta^{N_F^{gh}}.
		\end{align*}
		That is, $\{\widetilde{F}_\Delta\}_{(0,\bar{\Delta}]}$ provides order-$N_F^{gh}$ Hessian accuracy of $F$ at $x_0$.
		
		$\hfill\qed$
	\end{proof}
	
	Theorem \ref{thm:cha} shows that, in terms of Hessian accuracy at $x_0$, $\{\widetilde{F}_\Delta\}_{(0,\bar{\Delta}]}$ behaves at least as well as the worst behaved sub-model, where both gradient and Hessian behaviors are considered.
	
	\begin{theorem}{(Hessian accuracy near $x_0$)}\label{thm:chn}
		Let $f_1,f_2$ be in $\mathcal{C}^2$. Suppose $\{\widetilde{f}_{1,\Delta}\}_{(0,\bar{\Delta}]}$ provides order-$N_1^f$ function accuracy of $f_1$ near $x_0$, order-$N_1^g$ gradient accuracy of $f_1$ near $x_0$, and order-$N_1^h$ Hessian accuracy of $f_1$ near $x_0$. Suppose $\{\widetilde{f}_{2,\Delta}\}_{(0,\bar{\Delta}]}$ provides order-$N_2^f$ function accuracy of $f_2$ near $f_1(x_0)$, order-$N_2^g$ gradient accuracy of $f_2$ near $f_1(x_0)$, and order-$N_2^h$ Hessian accuracy of $f_2$ near $f_1(x_0)$. If $N_1^f > 0$, $\{\widetilde{f}_{2,\Delta}\}_{(0,\bar{\Delta}]}\subseteq\mathcal{C}^3$, and $\{||\nabla^3\widetilde{f}_{2,\Delta}||\}_{(0,\bar{\Delta}]}$ is uniformly bounded near $f_1(x_0)$ by $M_{\nabla^3\widetilde{f}_2}(f_1(x_0))$, then $\{\widetilde{F}_\Delta\}_{(0,\bar{\Delta}]}$ provides order-$N_F^{fgh}$ Hessian accuracy of $F$ near $x_0$, where $N_F^{fgh}=\min\{N_1^f,N_2^f,N_1^g,N_2^g,N_1^h,N_2^h\}$.
	\end{theorem}
	
	\begin{proof}
		Let $E_F^h = ||\nabla^2 F(x) - \nabla^2\widetilde{F}_\Delta(x)||$. Changing $x_0$ to $x$ in the proof of Theorem \ref{thm:cha}, we obtain
		\begin{align*}
		    E_F^h =& ||(\nabla\widetilde{f}_{2,\Delta}\big|_{f_1}^\top\nabla^2\widetilde{f}_{1,\Delta} - \nabla\widetilde{f}_{2,\Delta}\big|_{\widetilde{f}_{1,\Delta}}^\top\nabla^2\widetilde{f}_{1,\Delta} + \nabla f_1^\top\nabla^2\widetilde{f}_{2,\Delta}\big|_{f_1}\nabla f_1 - \nabla f_1^\top\nabla^2\widetilde{f}_{2,\Delta}\big|_{\widetilde{f}_{1,\Delta}}\nabla\widetilde{f}_{1,\Delta}\\
			&+ \nabla f_1^\top\nabla^2\widetilde{f}_{2,\Delta}\big|_{\widetilde{f}_{1,\Delta}}\nabla\widetilde{f}_{1,\Delta} - \nabla\widetilde{f}_{1,\Delta}^\top\nabla^2\widetilde{f}_{2,\Delta}\big|_{\widetilde{f}_{1,\Delta}}\nabla\widetilde{f}_{1,\Delta} + \nabla f_2\big|_{f_1}^\top\nabla^2\widetilde{f}_{1,\Delta} - \nabla\widetilde{f}_{2,\Delta}\big|_{f_1}^\top\nabla^2\widetilde{f}_{1,\Delta}\\
			&+ \nabla f_2\big|_{f_1}^\top\nabla^2 f_1 - \nabla f_2\big|_{f_1}^\top\nabla^2\widetilde{f}_{1,\Delta} + \nabla f_1^\top\nabla^2f_2\big|_{f_1}\nabla f_1 - \nabla f_1^\top\nabla^2\widetilde{f}_{2,\Delta}\big|_{f_1}\nabla f_1)(x)||.
		\end{align*}
		Noticing that 
		\begin{align*}
			|| \nabla f_2(f_1(x)) -  \nabla \widetilde{f}_{2,\Delta}(f_1(x))|| \le& \overline{M}_{\nabla f_1}(x_0)^{N_2^g}\kappa_2^g(f_1(x_0))\Delta^{N_2^g},\\
			|| \nabla^2 f_2(f_1(x)) -  \nabla^2 \widetilde{f}_{2,\Delta}(f_1(x))|| \le& \overline{M}_{\nabla f_1}(x_0)^{N_2^h}\kappa_2^h(f_1(x_0))\Delta^{N_2^h},
		\end{align*}
		and 
		\begin{align*}
			||\nabla\widetilde{f}_{2,\Delta}(f_1(x)) - \nabla \widetilde{f}_{2,\Delta}(\widetilde{f}_{1,\Delta}(x))|| \le& M_{\nabla^2\widetilde{f}_2}(f_1(x_0))\kappa_1^f(x_0)\Delta^{N_1^f},\\
			||\nabla^2\widetilde{f}_{2,\Delta}(f_1(x)) - \nabla^2 \widetilde{f}_{2,\Delta}(\widetilde{f}_{1,\Delta}(x))|| \le& M_{\nabla^3\widetilde{f}_2}(f_1(x_0))\kappa_1^f(x_0)\Delta^{N_1^f},
		\end{align*}
		we can get the conclusion.
		
		$\hfill\qed$
	\end{proof}
	
	Theorem \ref{thm:chn} shows that, in terms of Hessian accuracy near $x_0$, $\{\widetilde{F}_\Delta\}_{(0,\bar{\Delta}]}$ behaves at least as well as the worst behaved sub-model, where all function, gradient, and Hessian behaviors are considered.
	
	The next example shows the necessity of $N_1^f > 0$ in this section. In particular, if $N_1^f = 0$, then the function error of $\{\widetilde{F}_\Delta\}_{(0,\bar{\Delta}]}$ can be infinitely large.
	
	\begin{example}
		Suppose $\{\widetilde{f}_{1,\Delta}\}_{(0,\bar{\Delta}]}$ provides order-$N_1^f$ function accuracy of $f_1$ at $x_0$, and $\{\widetilde{f}_{2,\Delta}\}_{(0,\bar{\Delta}]}$ provides order-$N_2^f$ function accuracy of $f_2$ at $f_1(x_0)$. If $N_1^f = 0$, then $|F(x_0) - \widetilde{F}_\Delta(x_0)|$ can be infinitely large.
	\end{example}

	\begin{proof}
		Let $f_1(x)=1, \widetilde{f}_{1,\Delta}(x) = \Delta, f_2(x) = \frac{1}{x}, \widetilde{f}_{2,\Delta}(x) = \frac{1}{x}$. Then, regardless of $N_2^f$, we have
		\begin{align*}
			|f_2(x_0) - \widetilde{f}_{2,\Delta}(x_0)| =& 0 \le \Delta^{N_2^f},\\
			|F(x_0) - \widetilde{F}_\Delta(x_0)| =& |1 - \frac{1}{\Delta}|.
		\end{align*}
		Since $\lim\limits_{\Delta\to0} |1 - \frac{1}{\Delta}| =\infty$, we complete the proof.
		
		$\hfill\qed$
	\end{proof}

	In this example, $\{\widetilde{f}_{1,\Delta}\}_{(0,\bar{\Delta}]}$ and $\{\widetilde{f}_{2,\Delta}\}_{(0,\bar{\Delta}]}$ also obtain the true gradient and Hessian of $f_1$ and $f_2$. Therefore, even if we have more information in some theorems in this section, we still need $N_1^f > 0$.

	\section{Conclusion}
	
	In this paper, we examine the relation between the accuracy of sub-models and the accuracy of their product, quotient, and composition. This is analyzed in terms of order-$N$ terminology for function, gradient, and Hessian accuracy. Results show that the overall accuracy order is at least the same as the lowest accuracy order of sub-models. As corollaries, we demonstrate some situations where accuracy orders can be improved. We develop a broad framework that encompass many previous works and can be used as a theoretical tool in model error analysis. We also provide some examples to demonstrate the necessity of the assumption in our theorems. Notice that, throughout this paper, all global assumptions can be reduced to local assumptions.

    It should be mentioned that there are some cases that are not included in our framework. Our framework are based on the assumption that all functions are sufficiently smooth. In \cite{hare2013derivative,larson2016manifold}, the authors analyze sub-gradients and present error results for cases including nonsmooth component functions. In \cite{ling2021efficient,zhang2013structural}, an approximation method for performance function is developed. The authors use the product of several integrals to estimate fractional moments, which is one of the main steps in approximation. Although these cases are not included in this paper, they are possibly profitable follow-up research directions in our opinions.

    \section*{Acknowledgement}
    
    All authors' research is partially funded by the Natural Sciences and Engineering Research Council (NSERC) of Canada, Discover Grant \#2018-03865. Chen's research is partially funded by the MITACS Globalink program.
    \normalsize
    
    \bibliographystyle{siam}
    \bibliography{references}
    
\end{document}